\documentclass[letterpaper, 11pt]{article}

\usepackage{mathdots}
\usepackage{amsmath}
\usepackage{amssymb}
\usepackage{amsthm}
\usepackage{hyperref}
\usepackage{color}
\usepackage{caption}
\usepackage{calligra}
\usepackage{calrsfs}
\usepackage[mathscr]{euscript}

\usepackage[top=1in,left=1.5in,right=1.5in,bottom=1in]{geometry}

\newtheorem{theorem}{Theorem}[section]
\newtheorem{corollary}[theorem]{Corollary}
\newtheorem{lemma}[theorem]{Lemma}

\newtheorem{proposition}[theorem]{Proposition}
\theoremstyle{definition}
\newtheorem{definition}[theorem]{Definition}
\theoremstyle{remark}
\newtheorem{remark}[theorem]{Remark}
\newtheorem{example}[theorem]{Example}

\newcommand{\fop}{{\bf f}}
\newcommand{\gop}{{\bf g}}
\newcommand{\hop}{{\bf h}}

\newcommand{\cc}{{\bf c}}
\newcommand{\dd}{{\bf d}}

\title{On the Hurwitz Stability of Hurwitz-Type Matrix Polynomials 
}

\author{
Abdon E. Choque-Rivero \thanks{Email: \texttt{abdon.choque@umich.mx}}\\
Institute of Physics and Mathematics\\ Michoacan University of Saint Nicholas of Hidalgo\\
Morelia, Michoac\'an, M\'exico 
}

\date{}

\begin{document}

\maketitle

\begin{abstract}
Every matrix polynomial $\fop_n$ admits a decomposition of the form
\[
\fop_n(z)=\hop_n(z^2)+z\,\gop_n(z^2).
\]
The matrix polynomial $\fop_{2m}$ is said to be of Hurwitz type if the expression
$\gop_{2m}(z)\hop_{2m}^{-1}(z)$ admits a representation as a finite continued fraction with positive definite matrix coefficients. Similarly, the odd-degree matrix polynomial $\fop_{2m+1}$ is of Hurwitz type if
$\frac{1}{z}\hop_{2m+1}(z)\gop_{2m+1}^{-1}(z)$ has the same property.

We derive an explicit representation of the Bezoutian associated with Hurwitz-type matrix polynomials.
Using this representation, we obtain a direct proof that every Hurwitz-type matrix polynomial is Hurwitz.
The Hurwitz property of this class was also investigated in \cite{zhan1}; our approach is based on an explicit Bezoutian representation.
This provides a constructive connection between continued-fraction representations, matrix Bezoutians, and Hurwitz stability.
We also develop a completion procedure that associates with a given matrix polynomial a Hurwitz-type matrix polynomial of higher degree.
As a consequence, whenever such a completion exists, the original polynomial is Hurwitz.
The proposed construction is illustrated by examples.
\end{abstract}




\section{Introduction}

A scalar polynomial with all its roots in the left half of the
complex plane is called a Hurwitz polynomial
 \cite{routh, post, hurw,  datta, furhmann, holtz, VEM, ab2018, jaro, jury}. 
 The Hurwitz stability of such polynomials is also referred to as Hurwitzness---a term we adopt in this work.
 Let $n$ and $q$ be natural numbers.
 In the present work, we consider the properties and Hurwitzness of a subset of Lambda matrices \cite[page 95]{macd}, \cite[page 228]{gant01},
\cite{dennis}, defined by
\begin{equation} \fop_n(z):=A_{0} z^n+A_{1}
z^{n-1}+\ldots+  A_{n}, \label{hup01}
\end{equation}
 where $A_j$ is a complex $q\times q$ matrix and $z$ is a complex variable. 
As in \cite{lan01}, we refer to a matrix of the form (\ref{hup01}) as a matrix polynomial, which is
  called a Hurwitz matrix polynomial if $\det \fop_n$ is a Hurwitz polynomial.

The matrix polynomial (\ref{hup01})  has degree $n$ if $A_0$ is nonzero. Throughout this paper, we assume that $\det A_0\neq0$.
Clearly,  (\ref{hup01}) can be represented as
\begin{equation}
\fop_n(z)=   \hop_n(z^2)+z\,\gop_n(z^2), \label{hup02}
\end{equation}
where
  \begin{align}
 \hop_n(z)&:=
 \left\{
  \begin{array}{ll}
    A_{0} z^m+ A_{2} z^{m-1}+\ldots+ A_{2m},  & n=2m, \\
    A_{1}z^m+A_{3}z^{m-1}+\ldots+ A_{2m+1}, & n=2m+1,
  \end{array}
 \right. \label{hn}
 \end{align}
 and
 \begin{align}
\gop_n(z)&:=\left\{
  \begin{array}{ll}
    A_{1} z^{m-1}+ A_{3} z^{m-2}+\ldots+ A_{2m-1},  & n=2m, \\
    A_{0} z^m+A_{2}z^{m-1}+\ldots+ A_{2m}, & n=2m+1.
  \end{array}
 \right. \label{gn}
  \end{align}
 We will refer to
\eqref{hn} as the $\hop_n$ part and to \eqref{gn} as the $\gop_n$ part
of the matrix polynomial $\fop_n$.

Let us introduce the following notation. The symbol ${\mathbb C}^{q\times q}$ represents the set of all complex $q\times q$
matrices. For the null matrix that belongs to ${\mathbb C}^{p\times q}$, we will write $0_{p\times q}$. We denote by $0_q$ and $I_q$ the
null and identity matrices in ${\mathbb C}^{q\times q}$,
respectively. When the sizes of the null and identity matrices are
clear from the context, we will omit their indices.

In the following definition, we recall the concept of a Hurwitz-type matrix polynomial (HTM polynomial) introduced in \cite{abH}. In \cite{abH}, HTM polynomials were referred to as matrix Hurwitz-type polynomials.

For $A,B\in{\mathbb C}^{q\times q}$ with $B$ invertible, set
\begin{equation}\label{ABinv}
\frac{A}{B}:=AB^{-1}.
\end{equation}
\begin{definition}\label{defhur}
The $q\times q$ matrix polynomial $\fop_n$ in (\ref{hup01}) is called
a Hurwitz-type matrix polynomial if there exist two sequences of
positive definite matrices, $(\cc_{k})_{k=0}^{m-1}$ 
and
$(\dd_{k})_{k=0}^{m-1}$ (resp.  $(\cc_{k})_{k=0}^{m}$ and
$(\dd_{k})_{k=0}^{m-1}$) such that, for $n=2m$,
\begin{align}
 \frac{\gop_n (z)}{\hop_n(z)}=
 &\cfrac{I_q}{z\cc_{0}+\cfrac{I_q}{\dd_{0}
 +\cfrac{I_q}{\ddots 
 +\dd_{m-2}+
 \cfrac{I_q}{
+z\cc_{m-1}+\dd_{m-1}^{-1}}}}}, \label{cc13aa}
\end{align}
for all $z\in {\mathbb C}$ with $\det \hop_n(z)\neq0$ and, for $n=2m+1$,
\begin{align}
\frac{\hop_n (z)}{z \gop_n(z)}=
&\cfrac{I_q}{z\cc_{0}+\cfrac{I_q}{\dd_{0}+\cfrac{I_q}{\ddots \,
 +z\cc_{m-1}+\cfrac{I_q}{\dd_{m-1}
+z^{-1}\cc_{m}^{-1} }}}} \label{cc23aa}
 \end{align}
 for all $z\in {\mathbb C}\setminus\{0\}$ with $\det \gop_n(z)\neq0$.
\end{definition}
The definition of HTM polynomials in the form of a finite continued
fraction is inspired by  \cite[Chapter XV, Theorem
15]{gant0}.

\noindent
{\bf Main Results.} Below are the main results of the present work:
\begin{itemize}
\item[(a)]
We compute explicitly the Bezoutian of the HTM polynomial~(\ref{hup01}); see Theorem~\ref{th.3} and Corollary~\ref{cor3.5}.
From this, an explicit proof of the Hurwitz property of HTM polynomials readily follows; see Theorem~\ref{th.3aa}.
\item[(b)] For a given matrix polynomial $P_n$ of degree $n$ that is not a Hurwitz-type matrix (HTM) polynomial, we construct a polynomial
$\fop_{2n}(z)=P_n(z^2)+z Q_{n-1}(z^2)$
of degree $2n$, with a suitably chosen $Q_{n-1}$, such that $\fop_{2n}$ becomes an HTM polynomial.
Then, by applying the result from \cite{dyu2018}, we conclude that $P_n$ itself is a Hurwitz matrix polynomial; see Theorem~\ref{complet}.
\end{itemize} 

\medskip
\noindent
\textit{Relation between (a) and (b).} The explicit computation of the Bezoutian in part~(a) provides the structural foundation for the Hurwitz property of HTM polynomials. This result is then exploited in part~(b), where a general matrix polynomial is embedded into an HTM polynomial of higher degree. The Hurwitz property established in part~(a), together with the construction in~(b), allows us to transfer stability from the HTM setting back to the original polynomial $P_n$.

\medskip
\noindent

The Hurwitz property of matrix polynomials via Bezoutians was studied in \cite{heinig, lerer1, lerer2}. 
In \cite{zhan1}, the Hurwitz property of HTM polynomials was also considered. However, the proof of \cite[Theorem~3.5]{zhan1}, which serves as the basis for the analysis in that work, is not fully detailed. 
More precisely, several steps in the proof are left implicit; see Remark~\ref{rem5.7A}. In particular, the argument is developed only for the case $n=2m$, while the case $n=2m+1$ is described as ``similar,'' without explicit justification.

A distinctive feature of the odd-degree case $n=2m+1$ is that the associated block Hankel matrices $H_{1,m}$ and $H_{2,m-1}$ have different dimensions. 
In contrast to the even case, this asymmetry prevents a direct application of the standard construction. 
To overcome this difficulty, we introduce a modified polynomial $\widetilde \gop_{2m+1}$ (see \eqref{tg2m1}) of degree $m-1$, which replaces the original polynomial $\gop_{2m+1}$ (of degree $m$) in the construction of the symmetrizer.  
This approach leads to a consistent formulation of the matrix ${\mathcal W}_{2m+1}^{(0)}$ (see \eqref{B2mp1} and Remark~\ref{remult1}) and provides an explicit treatment of the odd-degree case, which is typically addressed only heuristically in the literature. 
This aspect is not  discussed in \cite{zhan1}.

\medskip

Additionally, in contrast to \cite{zhan1}, we use representations of the matrix polynomial 
$\fop_n$ in terms of orthogonal polynomials on $[0,\infty)$ of the first kind and polynomials of the second kind, as in (\ref{hup02}) and (\ref{hupfqp}), for example,
to compute explicitly the Bezoutian of the quadruple
$(\fop_n^*(\bar x), \fop_n(-x), \fop_n^*(-\bar x), \fop_n(x))$.

The asterisk denotes the Hermitian adjoint of a matrix, while the overbar represents complex conjugation.
In the scalar case, the use of the $\hop_n$ part
and $\gop_n$ part of the polynomial $\fop_n$ for studying the
corresponding Bezoutian was employed by Krein and Naimark in
\cite{kre-nai} and by Krein and Akhiezer in \cite{kre-akh}.

The importance of HTM polynomials as a subset of Hurwitz matrix polynomials is that they help study
  the asymptotic stability of the differential equation
$$
  I_q  y^{(n)}(t)+A_1 y^{(n-1)}(t)+\ldots+A_{n-1}y^\prime(t)+A_m y(t)=0, \ \ t\geq0
  $$
    with $y(t)\in {\mathbb R}^n$;  $I_q, A_1, \ldots, A_n\in {\mathbb C}^{q\times q}$; see
     \cite{martins, shieh, resende, henrion} and \cite{kraus}.

Our motivation is to generalize the results regarding  robust stability
\cite{choque2b} and finite-time stabilization with bounded controls
 \cite{choque2004, choque2004-2, choque5, choque6}  for the matrix case within the framework of HTM
polynomials.

 Our work is organized as follows. In Section 2, we review previous results on HTM polynomials and their associated Markov parameters. We also recall Stieltjes orthogonal matrix polynomials, orthogonal matrix polynomials, the interrelations between HTM polynomials and Markov parameters, as well as the Bezoutian and the inertia of matrix polynomials. 
In Section 3, we explicitly compute the Bezoutian of HTM polynomials. 
In Section 4, we study the Hurwitzness of HTM polynomials.
In Section 5, we present a result and an algorithm for constructing an HTM polynomial from a given matrix polynomial that is not of HTM type. We also explain why the original matrix polynomial is a Hurwitz matrix polynomial. 
In Section 6, we clarify certain steps left implicit in Theorem 3.5 of \cite{zhan1}, concerning a ``matrix analogue of the stability criterion via Markov parameters,'' as stated on page 10 of \cite{zhan1}.

This work is a direct continuation of the author's research in \cite{ablaa} and \cite{abH}, which focus on the nondegenerate truncated matricial Stieltjes problem and HTM polynomials, respectively.

\section{Preliminaries}\label{sec002}

In this section, we introduce the notation and preliminary results that will be used throughout the paper.
\subsection{HTM polynomials and associated Markov parameters}
The following quantities will be used to describe the domains of validity of the Laurent expansions appearing below.

If \(n=2m\), we define
\begin{equation}\label{nn01}
{\mathcal N}_{\hop_n}
:=
\{z\in {\mathbb C}: \det \hop_n(z)=0\},
\qquad
\eta_n
:=
\max\{|z|:\ z\in {\mathcal N}_{\hop_n}\},
\end{equation}
whereas, for \(n=2m+1\), we define
\begin{equation}\label{nn02}
{\mathcal N}_{\gop_n}
:=
\{z\in {\mathbb C}: \det \gop_n(z)=0\},
\qquad
\rho_n
:=
\max\{|z|:\ z\in {\mathcal N}_{\gop_n}\}.
\end{equation}

The next remark recalls the Laurent expansions of
\(
\gop_n(z)\hop_n(z)^{-1}
\)
and
\(
\hop_n(z)(z\gop_n(z))^{-1},
\)
which connect the matrix polynomial \eqref{hup01} with the associated Markov parameters \(s_j\).

\begin{remark}\label{rem2.1BB}
Let \(\hop_n\) and \(\gop_n\) be defined as in \eqref{hn} and \eqref{gn}, respectively, and let \(\fop_n\) be given by \eqref{hup02}.

For \(|z|>\eta_n\) (resp.\ \(|z|>\rho_n\)), the Laurent expansion in negative powers of \(z\) of
\(
\gop_n(z)\hop_n(z)^{-1}
\)
(resp.\ \(
\hop_n(z)(z\gop_n(z))^{-1}
\))
takes the form
\begin{align}
\frac{\gop_n(z)}{\hop_n(z)}
&=
\frac{s_0}{z}
-\frac{s_1}{z^2}
+\cdots
+(-1)^n\frac{s_n}{z^{n+1}}
+\cdots,
\qquad n=2m,
\label{hup09}
\\
\frac{\hop_n(z)}{z\gop_n(z)}
&=
\frac{s_0}{z}
-\frac{s_1}{z^2}
+\cdots
+(-1)^n\frac{s_n}{z^{n+1}}
+\cdots,
\qquad n=2m+1.
\label{hup10}
\end{align}
\end{remark}

For $ j\geq 0$ and by using the Markov parameters, we construct block Hankel matrices:
 \begin{equation} \label{61}
H_{1,j}:=\left(\begin{array}{cccc}
s_{0} & s_{1}& \ldots & s_{j}\\
s_{1} & s_{2}& \ldots & s_{j+1}\\
\vdots  & \vdots  & \vdots & \vdots\\
s_{j} & s_{j+1} & \ldots & s_{2j}
\end{array}\right),
H_{2,j}:=\left(\begin{array}{cccc}
s_{1} & s_{2}& \ldots & s_{j+1}\\
s_{2} & s_{3}& \ldots & s_{j+2}\\
\vdots  & \vdots  & \vdots & \vdots\\
s_{j+1} & s_{j+2} & \ldots & s_{2j+1}
\end{array}\right).
\end{equation}

The rational matrix functions in (\ref{cc13aa}) and (\ref{cc23aa}) are special solutions to the truncated Stieltjes matrix moment
(TSMM) problem; see proof of \cite[Theorem 7.9]{abH}. This implies that $\frac{\gop_n (-z)}{\hop_n(-z)}$  and
$\frac{\hop_n (-z)}{(-z) \gop_n(-z)}$ are Stieltjes transforms of a certain positive measure $\sigma$ on $[0,+\infty)$.

The TSMM problem is defined as follows: Let $(s_j)_{j=0}^{n}$ be a
finite sequence of $q\times q$ complex Hermitian matrices. The TSMM
problem consists of finding the set ${\mathcal M}_n$ of all
nonnegative Hermitian  $q\times q$ measures $\sigma$ defined on the
Borel $\sigma$--algebra ${\mathfrak B}\cap [0,+\infty)$  such that
$$
s_j=\int\limits_{[0,+\infty)} t^j \sigma(dt),\quad 0\leq j \leq n-1, \quad
s_n=\int\limits_{[0,+\infty)} t^n \sigma(dt)+M,
$$ 
 where $M$ is a $q\times q$ complex valued positive semi-definite matrix.

For each $\sigma \in {\mathcal M}_n$, one associates a matrix-valued function
\begin{equation}\label{szz}
s(z)=\int_{[0,+\infty)}\frac{\sigma(dt)}{t-z},
\end{equation}
which is defined and holomorphic in ${\mathbb C}\setminus [0,+\infty)$.

The function $s$ is associated with the TSMM problem and is called the Stieltjes transform of the positive measure $\sigma$; see
\cite{dyuviniti,dyu1999,dyu0,dyuth,dyu01,dyu2009,abmad,fmk01,fkm02,fkm03}.


In \cite{dyuviniti},  it is proven that  the TSMM problem in the case of an odd (resp.\ even) number of moments has a solution if and only if 
$H_{1,m}$ and $H_{2,m-1}$ (resp.\ $H_{1,m}$ and $H_{2,m}$) are both positive semidefinite matrices.

\begin{definition}\label{posseq}
The sequence $(s_{j})_{j=0}^{2m}$ (resp.\ $(s_{j})_{j=0}^{2m+1}$) is called  a Stieltjes positive sequence if both corresponding matrices
$H_{1,m}$ and $H_{2,m-1}$ (resp.\ $H_{1,m}$ and $H_{2,m}$) are positive definite.
In the case where $(s_{j})_{j=0}^{2m}$ (resp.\ $(s_{j})_{j=0}^{2m+1}$) is a positive Stieltjes sequence, the corresponding TSMM problem is called a TSMM problem in the nondegenerate case.
\end{definition}
In  \cite{dyu1999}, Yu. Dyukarev found the complete set of solutions to the TSMM problem in the case that $(s_{j})_{j=0}^{2m}$ and $(s_{j})_{j=0}^{2m+1}$ are Stieltjes positive sequences, respectively.

\subsection{Stieltjes orthogonal matrix polynomials} \label{sub1-2}
We now reproduce from \cite{ablaa} a number of auxiliary matrices, as well as two families of orthogonal
matrix polynomials and their polynomials of the second kind. Let
  $R_{j}:\mathbb{C}\to \mathbb{C}^{(j+1)q\times
(j+1)q}$
 be given by
$$
R_{j}(z):=(I_{(j+1)q}-zT_{j})^{-1},\qquad j\geq 0,
$$
 with
$$
T_{0}:=0_q, \ \  T_{j}:=\left(
\begin{array}{cc}
0_{q\times jq} & 0_q\\
I_{jq} & 0_{jq\times q}\\
                \end{array}
         \right),\qquad j\geq 1.
 $$
  Observe that for each $j\in {\mathbb N}_0$, the matrix--valued function
 $R_{j}$ can be represented via
\begin{equation}
R_{j}(z)=\left(
         \begin{array}{cccccc}
           I_q & 0_q & 0_q & \ldots & 0_q & 0_q \\
           zI_q &  I_q  & 0_q & \ldots & 0_q & 0_q\\
           z^2 I_q & z I_q & I_q & \ldots & 0_q & 0_q \\
 \vdots & \vdots & \vdots & \iddots & \vdots & \vdots \\
           z^jI_q & z^{j-1}I_q & z^{j-2}I_q &  \ldots& zI_q & I_q
         \end{array}
       \right).\label{o1.6}
\end{equation}
 Let
\begin{align}
v_{0}:=&I_q, 
\quad
v_{j}:=\left( \begin{array}{c}
I_q\\
0_{jq\times q}\\
\end{array}\right) \ \, \mbox{for all}\ \, j\in\mathbb{N}. \label{60a}
\end{align}
 Furthermore, let
$$
Y_{1,j}:=y_{[j,2j-1]},\,\  1\leq j\leq n, \,
                            \quad \mbox{and}
                            \quad Y_{2,j}:=y_{[j+1,2j]},\,\  2\leq j\leq
                            n,
$$
with
$$
y_{[j,k]}:=\left(s_j,s_{j+1},\ldots, s_k\right)^*,\quad 0\leq j\leq k 
$$
and  $y_{[j,k]}=0_q$ if $ j>k$.
 Define
\begin{align}
\label{66} &u_{1,0}:=0_q,\quad  u_{1,j}:=\left(\begin{array}{c}
0_q\\
-
y_{[0,j-1]}\\
\end{array}\right),\quad u_{2,j}:=\begin{array}{c}
-y_{[0,j]}\\
\end{array}, \quad 1\leq j\leq n.
\end{align}
Let $\widehat H_{1,j}$  (resp.\ $\widehat H_{2,j}$) denote the Schur complement of the block Hankel matrices $H_{1,j-1}$ in $H_{1,j}$ (resp.\,
 $H_{2,j-1}$ in $H_{2,j}$) as seen in \cite{ouellette}, \cite[Equalities (2.5) and (2.6)]{ablaa}:
\begin{align*}
   \widehat H_{1,0}:=s_0,\quad \widehat
   H_{1,j}:=&s_{2j}-Y_{1,j}^*H_{1,j-1}^{-1}Y_{1,j},\ \ j\geq1,
  \\
  \widehat H_{2,0}:=s_1,\quad \widehat H_{2,j}:=&s_{2j+1}-Y_{2,j}^*H_{2,j-1}^{-1}Y_{2,j},
  \ \  j\geq1.
 \end{align*}
 These matrices are   Hermitian positive definite matrices, and so
 are $H_{1,j}$ and $H_{2,j}$.

Now we partially reproduce left orthogonal matrix polynomials (OMP) on $[0,+\infty)$ \cite{dyu2018}.
\begin{definition}\label{def2pAA}
 Let $\sigma$ be a $q\times q$ nonnegenative measure on $[0,+\infty)$.
A sequence \((P_j)_{j=0}^\infty\) of \(q \times q\) matrix-valued polynomials is 
called orthogonal 
with respect to
 $\sigma$ 
 if 
\begin{item}
\item[i)] each matrix-valued polynomial \(P_j\) is of degree \(j\);
\item[ii)] the matrix-valued polynomials \(P_j\) satisfy the orthogonality relations:
\[
\int\limits_{[0,+\infty)} P_j(t) \, \sigma(dt) \, P_ k^*(t) = \delta_{jk} C_{qj}, 
\]
for all \(j, k \in \mathbb{N}\).
\end{item}
Here, \(\delta_{jk}\) is the Kronecker delta, and 
and $C_{qj}$ is a constant nonnegenative
 $q\times q$ matrix.
\end{definition}

The following matrix polynomials seem to have been first introduced in
\cite{dyuth}.
\begin{definition}
\label{def2p}
  Let $(s_{j})_{j=0}^{2n}$ (resp.\ $(s_{j})_{j=0}^{2n+1}$) be a Stieltjes positive sequence. Furthermore, let 
$v_j$, $u_{1,j}$, $u_{2,j}$ and $R_j$ be as in 
(\ref{60a}),  (\ref{66})  and (\ref{o1.6}), respectively. Let
$$
 P_{1,0}(z):=I_q, \ \ Q_{1,0}(z):=0_q, \ \ P_{2,0}(z):=I_q, \ \
Q_{2,0}(z):=s_0. 
$$
For $j\geq 1$,  set
\begin{align*}
P_{1,j}(z):=&(-Y_{1,j}^* H_{1,j-1}^{-1},\,I_q)R_j(z)v_j, \\
P_{2,j}(z):=&(-Y_{2,j}^*  H_{2,j-1}^{-1},\,I_q)R_j(z)v_j, \\
Q_{1,j}(z):=&-(-Y_{1,j}^*H_{1,j-1}^{-1},\,I_q)R_j(z)u_{1,j},\\
Q_{2,j}(z):=&-(-Y_{2,j}^*  H_{2,j-1}^{-1},\,I_q)R_j(z)u_{2,j}. 
\end{align*}
The matrix polynomials $Q_{1,j}$ and $Q_{2,j}$
 are called second-kind polynomials with respect to $P_{1,j}$ and  $P_{2,j}$,
respectively.
\end{definition}

In the next proposition, we recall that the matrix polynomials
$P_{1,j}$ and $P_{2,j}$ are orthogonal with respect to positive
measures supported on $[0,+\infty)$, and that the associated
second-kind polynomials $Q_{1,j}$ and $Q_{2,j}$ can be expressed
in terms of $P_{1,j}$ and $P_{2,j}$, respectively.
For completeness, we also record a basic identity involving
$P_{k,j}$ and $Q_{k,j}$.

 \begin{proposition} \label{p007}  a)
The  polynomials $P_{1,j}$ and $P_{2,j}$ are OMP with
 respect to $\sigma$ and $t\sigma$, respectively. More precisely (as seen in \cite[Remark 6.3]{abmad}),
$$
  \int\limits_{[0,+\infty)} P_{k,j}(t)t^{k-1}\sigma(dt)P_{k,l}^*(t)=\left\{\begin{array}{ccc}
                                                           0_q &: & j\neq l \\
                                                           \widehat
                                                           H_{k,j}&: &
                                                           j=l
                                                         \end{array}\right.,\ \  k=1,2.
  $$
 b)
 The following identities hold (as seen in \cite[Remark D.6]{ablaa} and \cite[Remark E.4]{ablaa}):
\begin{align}
 Q_{1,j}(x)=&\int\limits_{[0,+\infty)}
 \frac{P_{1,j}(x)-P_{1,j}(t)}{x-t}\sigma(dt), \ \  0\leq j\leq n,\label{qp1}\\
  Q_{2,j}(x)=&\int\limits_{[0,+\infty)}
  \frac{xP_{2,j}(x)-tP_{2,j}(t)}{x-t}\sigma(dt), \ \
  0\leq j\leq n-1.
 \nonumber 
   \end{align}
 c)  For $k=1,2$, the following identity holds (as seen in \cite[Proposition 5.1 (a)]{ablaa})
$$
P_{k,j}(z)Q_{k,j}^*(\bar z)-Q_{k,j}(z)P_{k,j}^*(\bar z)=0_q. 
$$
\end{proposition}
\vskip5mm

The following proposition establishes a connection between HTM polynomials and orthogonal matrix polynomials of the first and second kinds. The representation below follows from \cite[Theorem 6.1]{abH}, while identities \eqref{h001} and \eqref{g001} are established in the proof of \cite[Proposition 7.8]{abH}.

\begin{proposition}\label{prop-hupfqp}
For $n\ge 1$, every HTM polynomial admits the representation
\begin{equation}
\fop_n(z)=\left\{
\begin{array}{cl}
(-1)^m(P_{1,m}^*(-\bar z^2)-z\, Q_{1,m}^*(-\bar z^2)), & n=2m, \\[1mm]
(-1)^m(Q_{2,m}^*(-\bar z^2)+z\, P_{2,m}^*(-\bar z^2)), & n=2m+1.
\end{array}
\right.
\label{hupfqp}
\end{equation}

Moreover, the identities
\begin{align}
\hop_n(z)
&=
(-1)^m P_{1,m}^*(-\bar z),
\qquad
\gop_n(z)
=
(-1)^{m+1} Q_{1,m}^*(-\bar z),
\qquad n=2m,
\label{h001}
\end{align}
and
\begin{align}
\hop_n(z)
&=
(-1)^m Q_{2,m}^*(-\bar z),
\qquad
\gop_n(z)
&=
(-1)^m P_{2,m}^*(-\bar z),
\qquad n=2m+1
\label{g001}
\end{align}
hold.
\end{proposition}

The matrix polynomials $P_{k,j}$ and $Q_{k,j}$ also appear in the representation of the extremal solutions of the TSMM problem.

In the proof of \cite[Theorem 6.1]{abH}, the so--called extremal solutions of the TSMM problem are represented by
\begin{align}\label{GGG1}
G_{\sigma_{2m,\max}}(z)
=
-\frac{Q_{1,j}^*(\bar z)}{P_{1,j}^*(\bar z)}
=
\frac{\gop_{2m}(-z)}{\hop_{2m}(-z)}
\end{align}
and
\begin{align}\label{GGG2}
G_{\sigma_{2m,\min}}(z)
=
-\frac{Q_{2,j}^*(\bar z)}{z P_{2,j}^*(\bar z)}
=
\frac{\hop_{2m+1}(-z)}{(-z)\gop_{2m+1}(-z)}.
\end{align}

Note that the matrix function $G_{\sigma_{2m,\max}}$ is defined for
all $z\in {\mathbb C}\setminus [0,\infty)$ with $|z|>\eta_{2m}$,
whereas $G_{\sigma_{2m,\min}}$ is defined for
all $z\in {\mathbb C}\setminus [0,\infty)$ with $|z|>\rho_{2m+1}$.
Here, $\eta_{2m}$ and $\rho_{2m+1}$ are defined in (\ref{nn01}) and (\ref{nn02}), respectively.

\subsection{The Bezoutian and inertia for matrix polynomials}\label{sub1-6}
In this subsection,
 we reproduce some notions related to the inertia for a matrix polynomial \cite{lerer1,lerer2}.
 The  inertia for scalar polynomials was considered in \cite[page 462]{lan}.
\\
Let $L(\lambda)=L_0 \lambda^\ell+L_1\lambda^{\ell-1}+\ldots+\lambda L_{\ell-1}+L_\ell $ be
a $q\times q$ matrix polynomial with coefficients $L_k\in {\mathbb C}^{q\times q}$.
If $L_0\neq 0_q$, then ${\rm deg}(L)=\ell$ is  the degree of  polynomial  $L$. Furthermore,
if $L_0=I_q$, then $L$ is called monic.
 Let $M_1, L_1, M$ and  $L$ be polynomial matrices such that
 $M_1(\lambda)L_1(\lambda)-M(\lambda)L(\lambda)=0$.
 The  Bezoutian associated with the quadruple $(M_1, L_1,M, L)$ is the block matrix
 \begin{equation}\label{BM1M}
    B_{M_1,M}(L,L_1)=\left(\Gamma_{j,k}\right)_{j,k=1}^{m-1,\ell-1}.
\end{equation}
The block entries $\Gamma_{j,k}$ of the quadruple $(M_1, L_1,M, L)$
are determined by following equality:
 \begin{equation}\label{GML1}
 \sum_{j,k=0}^{m=1,\ell-1}\Gamma_{jk}\lambda^j\mu^k=\frac{1}{\lambda-\mu}\left(M_1(\lambda) L_1(\mu)-M(\lambda) L(\mu)
 \right).
 \end{equation}
The quantity $\ell$ (resp. $m$) denotes the maximal degree of
$L(\lambda)$, $L_1(\lambda)$ (resp. $M(\lambda)$, $M_1(\lambda)$).
 Assume that there exists $\lambda_0\in{\mathbb C}$ such that $\det L(\lambda_0)\neq0$:  on this case $L$ is called regular.

 Let $L$ be a regular polynomial of degree $\ell$. The triple
 \begin{equation}\label{inertia}
 \gamma(L)=\left(\gamma_{+}(L), \gamma_{-}(L), \gamma_{0}(L)\right)
 \end{equation}
 is called inertia with respect to the real axis, where $\gamma_{+}(L)$ (resp. $\gamma_{-}(L)$) is the number of zeros of
 $\det (L(\lambda))$ in the open upper (resp. open lower)
 half-plane, and $\gamma_0(L)$ is the number of zeros of  $\det (L(\lambda))$ in the real axis.
Now we reproduce \cite[Corollary 2.2]{lerer1}.

 \begin{corollary}\label{cor1.19}
     Let $L(\lambda)$ be a matrix polynomial with an invertible leading coefficient. The spectrum of $L(\lambda)$
     lies in the upper half-plane if and only if there exists a regular matrix polynomial $L_1(\lambda)$ that satisfies
$$
     L_1^*(\bar \lambda)     L_1(\lambda)=     L^*(\bar \lambda)  L(\lambda)
  $$ 
      such that the matrix 
      \begin{equation}\label{BLL}
            {\displaystyle \frac{1}{i}B_{L_1^*,L^*}(L,L_1)}
      \end{equation}
 is positive definite.
    \end{corollary}
    Our strategy for analyzing the Hurwitzness of HTM polynomials $\fop_n$ is based on the spectral symmetries of the quadruple
\[
(\fop_n^*(\bar x), \fop_n(-x), \fop_n^*(-\bar x), \fop_n(x)).
\]
The next lemma describes the spectral symmetries associated with this quadruple.

\begin{lemma}\label{lem:spectral_symmetries}
Let
$
F(z)=\sum_{k=0}^n A_k z^{n-k},
\, A_k\in\mathbb C^{q\times q},
$
be a regular matrix polynomial, and let
$
F^*(z):=\sum_{k=0}^n A_k^* z^{n-k}.
$
If
\[
\sigma(F):=\{z\in\mathbb C:\det F(z)=0\},
\]
then
$
\sigma(F(-z))=-\sigma(F),\,
$
$
\sigma(F^*(\bar z))=\overline{\sigma(F)},
$
and
$
\sigma(F^*(-\bar z))=-\overline{\sigma(F)}.
$
\end{lemma}
\begin{proof}
By definition,
$
\det F(-\lambda)=0
$
if and only if
$
-\lambda\in\sigma(F).
$
Hence
$
\sigma(F(-z))=-\sigma(F).
$
Next,
$
\det F^*(\bar z)
=
\overline{\det F(z)}.
$
Therefore,
$
\det F^*(\bar\lambda)=0
$
if and only if
$
\det F(\lambda)=0.
$
Consequently,
$
\sigma(F^*(\bar z))
=
\overline{\sigma(F)}.
$
Combining the previous two identities yields
$
\sigma(F^*(-\bar z))
=
-\overline{\sigma(F)}.
$
\end{proof}
\begin{remark}\label{rem:spectral_symmetries}
The four matrix polynomials
$
F(z)$,
$F(-z)$,
$F^*(\bar z)$
and
$F^*(-\bar z)$
correspond to natural symmetries of the spectrum of \(F\).
\\
Indeed, \(F(-z)\) reflects the spectrum through the origin,
\(F^*(\bar z)\) reflects the spectrum across the real axis, and
\(F^*(-\bar z)\) reflects the spectrum across the imaginary axis.
\\
Consequently, the Bezout form
\begin{equation}\label{FF1}
\frac{
F^*(\bar x)F(-y)-F^*(-\bar x)F(y)
}{x-y}
\end{equation}
compares the spectrum of \(F\) with its reflection across the imaginary axis.
\\
In the Hurwitz stable case,
$
\sigma(F)\subset\{z\in\mathbb C:\Re z<0\},
$
whereas
$
\sigma(F^*(-\bar z))
\subset
\{z\in\mathbb C:\Re z>0\}.
$
Hence these two spectra are separated by the imaginary axis. This separation is closely related to the definiteness properties of the associated Bezoutian.
\end{remark}

\subsection{The generalized Anderson--Jury Bezoutian matrix}\label{sub1-6A}

To compute the matrix \eqref{BLL} associated with an HTM polynomial, we use the generalized Bezoutian matrix and form introduced by Anderson and Jury in \cite{ander}, which are defined for a
$q\times q$ rational matrix function $W(z)$.

Suppose that $W(z)$ admits two matrix factorizations:
$$
W(z)=A^{-1}(z)B(z)=D(z)C^{-1}(z).
$$
The generalized Anderson--Jury Bezoutian form associated with the quadruple
$(A,B,C,D)$ is defined by
\[
\Gamma(x,y)
=
\frac{1}{x-y}\left(A(x)D(y)-B(x)C(y)\right).
\]
Using the notation \eqref{ABinv}, for HTM polynomials, the applicability of this construction follows from the fact that the $\hop_n$-part 
and the $\gop_n$-part of $\fop_n$ determine the rational matrix functions
\begin{equation}\label{sf2m}
s_{\fop_{2m}}(z)=\frac{\gop_{2m}(-z)}{\hop_{2m}(-z)}
\end{equation}
and
\begin{equation}\label{sf2mp1}
s_{\fop_{2m+1}}(z)=\frac{\hop_{2m+1}(-z)}{(-z)\gop_{2m+1}(-z)}.
\end{equation}
Taking into account that the functions \eqref{sf2m} and \eqref{sf2mp1} (see also \eqref{GGG1} and \eqref{GGG2}) are Stieltjes transforms of positive matrix-valued measures; see \eqref{szz}, they satisfy the Hermitian symmetry relation
\begin{equation}\label{sfn}
s_{\fop_n}(z)=s_{\fop_n}^*(\bar z).
\end{equation}
Therefore, the Anderson--Jury Bezoutian construction can be applied.
Consequently, rather than studying the Bezout form \eqref{FF1}, we employ Bezoutian forms associated with the quadruples
$
(x\gop_n^*, \hop_n, \hop_n^*, x\gop_n)
$
and
$
(\gop_n^*, \hop_n, \hop_n^*, \gop_n),
$
corresponding to the cases $n=2m+1$ and $n=2m$, respectively.
See the definitions of ${\mathcal F}_n$, ${\mathcal G}_n^{(1)}$, and ${\mathcal G}_n^{(2)}$ in \eqref{fn41}, \eqref{hgn1}, and \eqref{hgn2}, respectively, together with the relations between these forms given in \eqref{FGn}.

    \section{Bezoutian and Hurwitz-type  matrix  polynomials}\label{sec-4}

In this section, using the so-called matrix symmetrizers (\ref{Sqq}) and block Hankel matrices $H_{1,j}$ and $H_{2,j}$, we derive an explicit form of the Bezoutian for the quadruple
 ($\fop_n^*(\bar x),\fop_n(-x),\fop_n^*(-\bar x),\fop_n(x)$).

 Instead of working directly with the Bezoutian form
\begin{align}
{\mathcal F}_n(x,y):=&\frac{\fop_n^* (\bar x)\fop_n(-y)-\fop_n^*(-\bar x)\fop_n(y)}{x-y},\label{fn41}
\end{align}
we use the Equality (\ref{hup02}) and calculate the Bezoutian forms
\begin{align}
{\mathcal G}_n^{(1)}(x,y):=& \frac{1}{x^2-y^2}\left(
x^2 \gop_n^*(\bar x^2)\hop_n(y^2)- y^2 \hop_n^*(\bar x^2)\gop_n(y^2)\right)
\label{hgn1}
\end{align}
and
\begin{align}
{\mathcal G}_n^{(2)}(x,y):=& \frac{xy}{x^2-y^2} \left(
\gop_n^*(\bar x^2)\hop_n(y^2)-\hop_n^*(\bar x^2)\gop_n(y^2) \right) \label{hgn2}
\end{align}
associated with the quadruples $(x\gop_n^*, \hop_n, \hop_n^*, x\gop_n)$ and $(\gop_n^*, \hop_n, \hop_n^*, \gop_n)$.
Recall that a similar decomposition of the Bezoutian for scalar polynomials was previously studied; see \cite[page 293]{kre-nai} and \cite[page 257]{kre-akh}.
Following \cite{ander} to compute the Bezoutian forms ${\mathcal G}_n^{(1)}(x,y)$ and ${\mathcal G}_n^{(2)}(x,y)$
we
employ both representations (\ref{sfn}) of the rational matrices (\ref{cc13aa}) and (\ref{cc23aa}).


It readily  follows that
\begin{equation}\label{FGn}
{\mathcal F}_n(x,y)=2\left({\mathcal G}_n^{(1)}(x,y)+{\mathcal G}_n^{(2)}(x,y)\right).
\end{equation}

In the next lemma, we demonstrate that the matrix functions
${\mathcal F}_n$, ${\mathcal G}_n^{(1)}$ and ${\mathcal G}_n^{(2)}$ are matrix polynomials in $x$ and $y$. See \cite[Remark 1]{ander}.

\begin{lemma}\label{lem3.1A}
Let ${\mathcal F}_n$,
 ${\mathcal G}_n^{(1)}$ and ${\mathcal G}_n^{(2)}$ be as defined in (\ref{fn41}), (\ref{hgn1}) and (\ref{hgn2}),
 respectively.
 The functions  ${\mathcal F}_n$,  ${\mathcal G}_n^{(1)}$ and ${\mathcal G}_n^{(2)}$ are polynomials on $x$ and $y$.
\end{lemma}
\begin{proof} 
Using the numerator on the right-hand side of \eqref{fn41}, together with \eqref{hupfqp}, and replacing $y$ by $x$, we obtain
\begin{align*}
&\fop_n(\bar x)\fop_n(-x)-\fop_n(-\bar x)\fop_n(x)\\
&=2x\left\{\begin{array}{lc}
P_{1,m}(-x^2)Q_{1,m}^*(-x^2)-Q_{1,m}(-x^2)P_{1,m}^*(-x^2),& n=2m\\
P_{2,m}(-x^2)Q_{2,m}^*(-x^2)-Q_{2,m}(-x^2)P_{2,m}^*(-x^2),& n=2m+1
\end{array}
\right.\\
&=0.
\end{align*}
By \cite[Chapter~1]{CLO}, any polynomial vanishing on the diagonal $x=y$
is divisible by $x-y$. Since $\fop_n^* (\bar x)\fop_n(-y)-\fop_n^*(-\bar x)\fop_n(y)$ has this property,  ${\mathcal F}_n(x,y)$ is a polynomial.
Arguing similarly with the numerators in \eqref{hgn1} and \eqref{hgn2}, and following \cite[Chapter~1]{CLO}, 
we show that ${\mathcal G}_n^{(1)}$ and ${\mathcal G}_n^{(2)}$ are polynomials in $x$ and $y$.
\end{proof}
   
 For $k=0,1$, 
 let
 \begin{align}
H^{(k)}_\infty:=& \begin{pmatrix}
s_k&-s_{k+1}&s_{k+2}&\ldots\\
-s_{k+1}&s_{k+2}&\iddots\\
s_{k+2}&\iddots
\end{pmatrix}
  \label{Hk0}
\end{align}
be an infinite Hankel matrix.

A similar untruncated block Hankel matrix, denoted by $\bf {H}$,
appears in the proof of  \cite[Lemma 2.3]{ander}.
\begin{lemma}\label{lem3.2}
Let ${\mathcal F}_n$, ${\mathcal G}_n^{(1)}$ and ${\mathcal G}_n^{(2)}$ be as in (\ref{fn41}),  (\ref{hgn1}) and
(\ref{hgn2}), respectively.
Thus, for $n=2m$, we have
\begin{align}
{\mathcal G}_{2m}^{(1)}(x,y)=&\hop_{2m}^*(\bar x^2)\left(x^{-2}I_q,x^{-4}I_q,\ldots\right)
H^{(1)}_\infty \left(\bar y^{-2}I_q,\bar y^{-4}I_q,\ldots\right)^*\hop_{2m}(y^2),
\label{G2m1}
\\
{\mathcal G}_{2m}^{(2)}(x,y)=&-\hop_{2m}^*(\bar x^2)\left(x^{-1}I_q,x^{-3}I_q,\ldots\right)
H^{(0)}_\infty \left(\bar y^{-1}I_q,\bar y^{-3}I_q,\ldots\right)^* \hop_{2m}(y^2).
\label{G2m2}
\end{align}
For $n=2m+1$, we have
\begin{align}
{\mathcal G}_{2m+1}^{(1)}(x,y)=&\gop_{2m+1}^*(\bar x^2)\left(I_q,x^{-2}I_q,\ldots\right)
H^{(0)}_\infty \left(I_q,\bar y^{-2}I_q,\ldots\right)^* \gop_{2m+1}(y^2),
\label{G2mp1}\\
{\mathcal G}_{2m+1}^{(2)}(x,y)=&-\gop_{2m+1}^*(\bar x^2)\left(x^{-1}I_q,x^{-3}I_q,\ldots\right)
H^{(1)}_\infty
\left(\bar y^{-1}I_q, \bar y^{-3}I_q,\ldots\right)^* \nonumber\\
&\cdot \gop_{2m+1}(y^2).
\label{G2mp2}
\end{align}
\end{lemma}
\begin{proof}
We prove (\ref{G2m1}). With (\ref{hup09}), we have
\begin{align*}
{\mathcal G}_{2m}^{(1)}(x,y)=&\frac{1}{x^2-y^2}\hop_{2m}^*(\bar x^2)\left(x^2\hop_{2m}^{-1^*}(\bar x^2)\gop_{2m}^*(\bar x^2)-y^2 \gop_{2m}(y^2)\hop_{2m}^{-1}(y^2)\right)\nonumber\\
&\cdot \hop_{2m}(y^2)\\
=&\frac{1}{x^2-y^2}\hop_{2m}^*(\bar x^2)\left(s_1(y^{-2}-x^{-2})-s_2(y^{-4}-x^{-4})+s_3(y^{-6}-x^{-6})\right.\nonumber\\
&\left.-\ldots\right)\hop_{2m}(y^2)\\
=&\frac{1}{x^2y^2}\hop_{2m}^*(\bar x^2)\left(s_1-s_2(x^{-2}+y^{-2})+s_3(x^{-4}+x^{-2}y^{-2}+y^{-4})\right.\nonumber\\
&\left.-\ldots\right)\hop_{2m}(y^2)\\
=&
\hop_{2m}^*(\bar x^2)
\left(x^{-2}I_q,x^{-4}I_q,\ldots\right)
H^{(1)}_\infty \left(\bar y^{-2}I_q, \bar y^{-4}I_q,\ldots\right)^*
 \hop_{2m}(y^2).
\end{align*}
In the last equality, we substituted (\ref{Hk0}) for $k=1$.
Similarly, to prove (\ref{G2m2}), we use (\ref{hup09}) and  (\ref{Hk0}) for $k=0$.
Equalities (\ref{G2mp1}) and  (\ref{G2mp2}) are proven by applying (\ref{hup10}) and (\ref{Hk0}) for $k=0$ and $k=1$, respectively.
\end{proof}
\begin{remark}\label{rem03}
The infinite block Hankel matrix $H^{(1)}_\infty$ appearing in \eqref{G2m1} arises naturally from the Laurent expansions in \eqref{hup09}. 
Indeed, the expression
\[
\frac{1}{x^2-y^2}
\left(x^2\hop_{2m}^{-1^*}(\bar x^2)\gop_{2m}^*(\bar x^2)-y^2 \gop_{2m}(y^2)\hop_{2m}^{-1}(y^2)\right)
\]
generates an infinite series in negative powers of $x$ and $y$, whose coefficients are precisely the moments $(s_j)_{j\geq 1}$, giving rise to the infinite block Hankel matrix $H^{(1)}_\infty$ in the representation
\[
\left(x^{-2}I_q,x^{-4}I_q,\ldots\right)
H^{(1)}_\infty \left(\bar y^{-2}I_q, \bar y^{-4}I_q,\ldots\right)^*.
\]
However, since ${\mathcal G}_{2m}^{(1)}$ is a polynomial, only finitely many terms contribute to the products appearing in \eqref{G2m1}. 
Consequently, the right-hand side of \eqref{G2m1} effectively reduces to a finite truncation of $H^{(1)}_\infty$.
\\
The same observation applies to the right-hand sides of
\eqref{G2m2}--\eqref{G2mp2}.
\end{remark}

\medskip
\noindent
In view of Remark~\ref{rem03}, the infinite block vector
$$
\left(\bar \xi^{-2}I_q,\bar \xi^{-4}I_q,\ldots\right)^*\hop_{2m}(\xi^2)
$$
can be decomposed into a polynomial part and a part containing strictly negative powers of $\xi$, namely,
\begin{equation}
\left(\bar \xi^{-2}I_q,\bar \xi^{-4}I_q,\ldots\right)^*\hop_{2m}(\xi^2)
=
\hop_{2m}^{(11)}(\xi)+ \hop_{2m}^{(12)}(\xi).
\label{hh11}
\end{equation}
Similarly, we have
\begin{align}
\left(\bar \xi^{-1}I_q,\bar \xi^{-3}I_q,\ldots\right)^*\hop_{2m}(\xi^2)
&=
\hop_{2m}^{(21)}(\xi)+ \hop_{2m}^{(22)}(\xi),
\label{hh21}
\\
\left(I_q,\bar \xi^{-2}I_q,\ldots\right)^* \gop_{2m+1}(\xi^2)
&=
\gop_{2m+1}^{(11)}(\xi)+ \gop_{2m+1}^{(12)}(\xi),
\label{gg11}
\\
\left(\bar \xi^{-1}I_q,\bar \xi^{-3}I_q,\ldots\right)^*\gop_{2m+1}(\xi^2)
&=
\gop_{2m+1}^{(21)}(\xi)+ \gop_{2m+1}^{(22)}(\xi).
\label{gg21}
\end{align}
Here, $\hop_{2m}^{(k1)}(\xi)$ (resp.\ $\gop_{2m+1}^{(k1)}(\xi)$), for $k=1,2$, are polynomials in $\xi$, whereas
$\hop_{2m}^{(k2)}(\xi)$ (resp.\ $\gop_{2m+1}^{(k2)}(\xi)$), for $k=1,2$, are polynomials in $\xi^{-1}$.

\medskip

For a natural number $m$, let
\begin{equation}\label{Jqq}
\widetilde J_m
:=
{\rm diag}\,\left(I_q,(-1)I_q,(-1)^2I_q,\ldots,(-1)^mI_q\right)
\end{equation}
be a $(m+1)q\times (m+1)q$ matrix.

The following remark explains how the infinite Hankel matrices $H^{(k)}_\infty$ reduce to finite block Hankel matrices when acting on the polynomial components
$\hop_{2m}^{(k1)}(\xi)$ (resp.\ $\gop_{2m+1}^{(k1)}(\xi)$), for $k=1,2$, introduced in \eqref{hh11}--\eqref{gg21}. This follows from direct calculations based on the corresponding decompositions.
\begin{remark}
Let
$\hop_{2m}^{(11)}(\xi)$, $\hop_{2m}^{(21)}(\xi)$, $ \gop_{2m+1}^{(11)}(\xi)$ and $ \gop_{2m+1}^{(21)}(\xi)$
be as defined in (\ref{hh11}), (\ref{hh21}), (\ref{gg11}) and (\ref{gg21}), respectively.
Furthermore, let $H_{1,j}$, $H_{2,j}$ and $J_m^{(k)}$ be as defined in (\ref{61}) and (\ref{Jqq}), respectively.
Thus, the following equalities hold
\begin{align}
\hop_{2m}^{(11)^*}(\bar x)H^{(1)}_\infty\hop_{2m}^{(11)}(y)=&
\widehat \hop_{2m}^{(11)^*}(\bar x)\widetilde J_{m-1} H_{2,m-1}\widetilde J_{m-1}\widehat \hop_{2m}^{(11)}(y), \label{hN1}
\\
\hop_{2m}^{(21)^*}(\bar x)H^{(1)}_\infty\hop_{2m}^{(21)}(y)=&
\widehat \hop_{2m}^{(21)^*}(\bar x)\widetilde J_{m-1}H_{1,m-1}\widetilde J_{m-1}\widehat \hop_{2m}^{(21)}(y), \label{hN2}
\\
\gop_{2m+1}^{(11)^*}(\bar x)H^{(1)}_\infty\gop_{2m+1}^{(11)}(y)=&
\widehat \gop_{2m+1}^{(11)^*}(\bar x)\widetilde J_{m}H_{1,m}\widetilde J_{m}\,\widehat \gop_{2m+1}^{(11)}(y), \label{hN3}
\\
\gop_{2m+1}^{(21)^*}(\bar x)H^{(1)}_\infty\gop_{2m+1}^{(21)}(y)=&
\widehat \gop_{2m+1}^{(21)^*}(\bar x)\widetilde J_{m-1}H_{2,m-1}\widetilde J_{m-1}\widehat \gop_{2m+1}^{(21)}(y), \label{hN4}
\end{align}
 where
$$
 \widehat \hop_{2m}^{(11)},\, \widehat \hop_{2m}^{(21)},\,  \widehat \gop_{2m+1}^{(11)} \quad \mbox{and} \quad \widehat \gop_{2m+1}^{(21)}
 $$
 represent finite-dimensional restrictions of the matrices
 $\hop_{2m}^{(11)}$, $\hop_{2m}^{(21)}$, $\gop_{2m+1}^{(11)}$ and $\gop_{2m+1}^{(21)}$.
 In particular, the polynomial part of
 \begin{equation}\label{hNN1}
 \widehat \hop_{2m}^{(11)}(y)=\left(\bar y^{-2}\hop_{2m}^*(\bar y^2), \bar y^{-4}\hop_{2m}^*(\bar y^2),\ldots,\bar y^{-2m}\hop_{2m}^*(\bar y^2)\right)^*
  \end{equation}
 is the restriction of the polynomial part of the matrix
 $$
 \hop_{2m}^{(11)}(y)=\left(\bar y^{-2}\hop_{2m}^*(\bar y^2), \bar y^{-4}\hop_{2m}^*(\bar y^2),\ldots, \bar y^{-2m}\hop_{2m}^*(\bar y^2),0_q,0_q,\ldots\right)^*.
 $$
\end{remark}
  Equalities (\ref{hN1})--(\ref{hN4}) correspond to those in the proof of \cite[Lemma 2.3]{ander}, where, instead of the matrices \( H_{1,j} \) and \( H_{2,j} \), the truncated block Hankel matrix \( {\bf H}_{nm} \) is used. Moreover, in place of the matrix \( H_\infty^{(1)} \), the infinite Hankel matrix \( H \) is employed.
\vskip2mm
For a matrix polynomial $Q(x)=Q_0 x^\ell+Q_1 x^{\ell-1}+\ldots+Q_{\ell-1}x+Q_\ell$, the Hankel block matrix
\begin{equation}\label{Sqq}
  S(Q):=\begin{pmatrix}
Q_{\ell-1}& Q_{\ell-2}&\ldots&Q_1&Q_0\\
Q_{\ell-2}& Q_{\ell-3}&\ldots&Q_0&0_q\\
  \vdots&  \iddots&  \iddots&  \iddots&  \vdots\\
  Q_0& 0_q&\ldots&0_q&0_q
 \end{pmatrix}
\end{equation}
is called a matrix symmetrizer of $Q(x)$. 
This matrix provides a finite-dimensional representation of the coefficients of $Q$ and will be used to express the Bezoutian in terms of block Hankel matrices.
In \cite[page 455]{lan}, the concept of symmetrizer is introduced for scalar polynomials.

Additionally, let
\begin{align}
 F_0^{2k}(x):=&\left(I_q, \bar x^2I_q,\ldots,\bar x^{2k}I_q\right)^*,\label{F02k}\\
  F_1^{2k-1}(x):=&\left(\bar xI_q,\bar x^3I_q,\ldots,\bar x^{2k-1}I_q\right)^*.\label{F12km1}
\end{align}
In the following theorem,  added by symmetrizers of corresponding matrices and the block Hankel matrices $H_{k,j}$, we 
obtain a factorized expression of the Bezoutian forms
${\mathcal G}_{2m}^{(1)}$,
${\mathcal G}_{2m}^{(2)}$,
${\mathcal G}_{2m+1}^{(1)}$ and
${\mathcal G}_{2m+1}^{(2)}$. 
%
\begin{theorem}\label{th.3}
Let ${\mathcal G}_{2m}^{(k)}$, ${\mathcal G}_{2m+1}^{(k)}$ for $k=1,2$ be as defined in (\ref{hgn1}) and (\ref{hgn2}).
Additionally, let  $\hop_{n}$,  $\gop_{n}$, $S$, $\widetilde J_m$, $H_{1,j}$, $H_{2,j}$,  $F_0^{2k}(x)$ and $  F_1^{2k-1}(x)$ be as in (\ref{hn}),  (\ref{gn}), (\ref{Sqq}),  (\ref{Jqq}),
(\ref{61}), (\ref{F02k}) and  (\ref{F12km1}),  respectively.\\
 Therefore,
for $n=2m$, the following equality holds:
\begin{align}
{\mathcal G}_{2m}^{(1)}(x,y)=&F_0^{(2m-2)^*}(\bar x)S(\hop_{2m}^*)\widetilde J_{m-1}  H_{2,m-1} \widetilde J_{m-1}  S(\hop_{2m})
F_0^{(2m-2)}(y),\label{GGa2m1}
\\
{\mathcal G}_{2m}^{(2)}(x,y)=&-F_1^{(2m-1)^*}(\bar x)S(\hop_{2m}^*)\widetilde J_{m-1} H_{1,m-1}\widetilde J_{m-1} S(\hop_{2m})F_1^{(2m-1)}(y).
\label{GGa2m2}
\end{align}
For $n=2m+1$, we have
\begin{align}
{\mathcal G}_{2m+1}^{(1)}(x,y)=&F_0^{(2m)^*}(\bar x)S(\gop_{2m+1}^*)\widetilde J_{m} H_{1,m}\widetilde J_{m}  S(\gop_{2m+1})
F_0^{(2m)}(y). \label{GGa2mp1}
\end{align}
Moreover,
\begin{align}
{\mathcal G}_{2m+1}^{(2)}(x,y)=&-F_1^{(2m-1)^*}(\bar x)S(\widetilde \gop_{2m+1}^*)\widetilde J_{m-1} H_{2,m-1}\widetilde J_{m-1} S(\widetilde \gop_{2m+1})
F_1^{(2m-1)}(y),\label{GGa2mp2}
\end{align}
where
\begin{equation}\label{tg2m1}
\widetilde \gop_{2m+1}(z):= I_q z^{m-1}+A_{2}z^{m-2}+\ldots+ A_{2m-2}.
\end{equation}
\end{theorem}
\begin{proof}
We prove (\ref{GGa2m1}). From (\ref{G2m1}), (\ref{hh11}) and  (\ref{hN1}), we have
\begin{align}
{\mathcal G}_{2m}^{(1)}(x,y)=\widehat \hop_{2m}^{(11)^*}(\bar x)\widetilde J_{m-1} H_{2,m-1}\widetilde J_{m-1} \widehat \hop_{2m}^{(11)}(y).
\label{hNa}
\end{align}
We examine $\widehat \hop_{2m}^{(11)}(y)$. By using (\ref{hNN1}) and (\ref{Sqq}), we have
\begin{align}
\widehat \hop_{2m}^{(11)}(y)=&
\begin{pmatrix}
 A_{2m-2}+\ldots+A_{2}y^{2m-4}+I_qy^{2m-2}\\
A_{2m-4}+\ldots+A_{2}y^{2m-6}+I_q y^{2m-4}\\
\ldots\\
I_q
\end{pmatrix}\nonumber\\
=&
\begin{pmatrix}
 A_{2m-2}&\ldots&A_{2}&I_q\\
A_{2m-4}&\ldots&I_q&\\
\ldots\\
I_q
\end{pmatrix}
\begin{pmatrix}
I_q\\
y^{2}I_q\\
\ldots\\
y^{2m-2}I_q
\end{pmatrix}\nonumber\\
=&S(\hop_{2m})\left(I_q,\bar y^2I_q,\ldots, \bar y^{2m-2}I_q\right)^*.\label{ySh2}
\end{align}
By substituting (\ref{ySh2}) for (\ref{hNa}), we obtain (\ref{GGa2m1}).
Similarly, by using 
\begin{align}
{\mathcal G}_{2m}^{(2)}(x,y)=&\widehat \hop_{2m}^{(21)^*}(\bar x)\widetilde J_{m-1}H_{1,m-1}\widetilde J_{m-1} \widehat \hop_{2m}^{(21)}(y),\label{hNb}\\
{\mathcal G}_{2m+1}^{(1)}(x,y)=&\widehat \gop_{2m+1}^{(11)^*}(\bar x)\widetilde J_{m} H_{1,m-1}\widetilde J_{m}\, \widehat \gop_{2m+1}^{(11)}(y),\label{hNc}\\
{\mathcal G}_{2m+1}^{(2)}(x,y)=&\widehat \gop_{2m+1}^{(21)^*}(\bar x)\widetilde J_{m-1} H_{2,m-1}\widetilde J_{m-1} \widehat \gop_{2m+1}^{(21)}(y),\label{hNd}
\end{align}
one can prove (\ref{GGa2m2}), (\ref{GGa2mp1}) and (\ref{GGa2mp2}).
Equalities (\ref{hNb}), (\ref{hNc}) and (\ref{hNd}) readily follow from (\ref{hN2}), (\ref{hN3}), (\ref{hN4}) and (\ref{hh21}), as well as from  (\ref{gg11}) and (\ref{gg21}).
Additionally, one also uses (\ref{G2m2}), (\ref{G2mp1}), (\ref{G2mp2}),  (\ref{F02k}) and (\ref{F12km1}).
\end{proof}
\begin{remark}\label{remult1}
The polynomial $\widetilde \gop_{2m+1}$ defined in \eqref{tg2m1} has degree $m-1$, in contrast to the original polynomial $\gop_{2m+1}$ defined in \eqref{gn}.
Using $\widetilde \gop_{2m+1}$, we construct the symmetrizer $S(\widetilde \gop_{2m+1})$ appearing in the matrix ${\mathcal W}_{2m+1}^{(0)}$ in \eqref{B2mp1}.
This is consistent with the fact that the block Hankel matrices $H_{1,m}$ and $H_{2,m-1}$ have different dimensions.
\end{remark}

In the following corollary, we represent ${\mathcal F}_n$ as a matrix quadratic form
 regarding the
  block Hankel matrices $H_{1,j}$ and $H_{2,j}$, as well as the symmetrizer matrices $S(\hop_{2m})$, $S(\gop_{2m+1})$ and $S(\widetilde \gop_{2m+1})$.
\begin{corollary}\label{cor3.5}
Assume the same conditions and notations as in Theorem \ref{th.3}.  Let $\widetilde  J_{m}$ be as in (\ref{Jqq}). 
Moreover, let
\begin{equation}\label{Jkl}
{\mathcal J}_{r}^{(\ell)}:=\begin{pmatrix}
\widetilde  J_{\ell}&0\\
0&\widetilde  J_{r}
 \end{pmatrix}
\end{equation}
for all $r,\,\ell>1$.
  Furthermore, for $k=0,1$ let
 \begin{align}
{\mathcal W}_{2m}^{(k)}:=&
\begin{pmatrix}
 S(\hop_{2m}^*)&0\\
0& S(\hop_{2m}^*)
 \end{pmatrix}
 {\mathcal J}_{m-1}^{(m-1)}
\begin{pmatrix}
 H_{2,m-1}&0\\
0& (-1)^k H_{1,m-1}
 \end{pmatrix}
  {\mathcal J}_{m-1}^{(m-1)}
 \nonumber\\
 &\cdot
  \begin{pmatrix}
 S(\hop_{2m})&0\\
0& S(\hop_{2m})
 \end{pmatrix}, \label{WW1}
 \\
{\mathcal W}_{2m+1}^{(k)}:=&
  \begin{pmatrix}
 S(\gop_{2m+1}^*)&0\\
0& S(\widetilde \gop_{2m+1}^*)
 \end{pmatrix}
  {\mathcal J}_{m-1}^{(m)}
 \begin{pmatrix}
 H_{1,m}&0\\
0& (-1)^k H_{2,m-1}
 \end{pmatrix}
   {\mathcal J}_{m-1}^{(m)}
 \nonumber
 \\
 &\cdot
  \begin{pmatrix}
 S(\gop_{2m+1})&0\\
0& S(\widetilde \gop_{2m+1})
 \end{pmatrix}.\label{WW2}
\end{align}
With ${\mathcal F}_n$ be as in (\ref{fn41}),
the following equalities hold:
\begin{align}
 {\mathcal F}_{2m}(x,y)=&2
 \begin{pmatrix}
 F_0^{(2m-2)}(\bar x)\\
 F_1^{(2m-1)}(\bar x)
 \end{pmatrix}^{*}
  {\mathcal W}_{2m}^{(1)}
 \begin{pmatrix}
 F_0^{(2m-2)}(y)\\
 F_1^{(2m-1)}(y)
 \end{pmatrix},\label{F2mC}
 \end{align}
 and
\begin{align}
 {\mathcal F}_{2m+1}(x,y)=&2
 \begin{pmatrix}
 F_0^{(2m)}(\bar x)\\
 F_1^{(2m-1)}(\bar x)
 \end{pmatrix}^{*}
  {\mathcal W}_{2m+1}^{(1)}
 \begin{pmatrix}
 F_0^{(2m)}(y)\\
 F_1^{(2m-1)}(y)
 \end{pmatrix}.\label{F2mp1C}
 \end{align}
\end{corollary}
\begin{proof}
Equality \eqref{F2mC} follows directly from
\eqref{FGn}, \eqref{GGa2m1}, \eqref{GGa2m2}, and \eqref{WW1}.
In a similar manner, Equality (\ref{F2mp1C}) can readily be derived from (\ref{FGn}), (\ref{GGa2mp1}), (\ref{GGa2mp2}) and (\ref{WW2}).
\end{proof}
%


    \section{Hurwitness of Hurwitz-type  matrix  polynomials}\label{sec05}

  In this section, using the explicit form of the Bezoutian matrix (\ref{F2mC}) and (\ref{F2mp1C}), we verify that an HTM polynomial $\fop_n$ is indeed a Hurwitz matrix polynomial.

In \cite{zhan1}, the Hurwitz property of HTM polynomials was studied. 
However, the proof of \cite[Theorem~3.5]{zhan1} concerning the inertia of the matrix polynomial $\fop_n$ for $n=2m$ is neither explicit nor complete, and a proof for the case $n=2m+1$ is not provided. See Remark~\ref{rem5.7A}.

     Recall that the inertia of matrix polynomials, denoted by  $\gamma(\fop_n)$, provides the number of eigenvalues of $\fop_n$
    with positive, negative and real zero parts, respectively; see \cite[page 410]{lerer1}.
\vskip2mm
We rewrite Corollary \ref{cor3.5} for the pair $(ix,iy)$ instead of $(x,y)$.
\begin{corollary}\label{cor5.1}
Assume the same conditions and notations as in Theorem \ref{th.3}.
 Moreover, %
 let ${\mathcal F}_n$,  ${\mathcal J}_k^\ell$, ${\mathcal W}_{2m}^{(0)}$ and
${\mathcal W}_{2m+1}^{(0)}$ be
as in (\ref{fn41}), (\ref{Jkl}), (\ref{WW1}) and (\ref{WW2}), respectively. Thus, the following equalities hold:
\begin{align}
 {\mathcal F}_{2m}(ix,iy)=&2
 \begin{pmatrix}
 F_0^{(2m-2)}(\bar x)\\
 F_1^{(2m-1)}(\bar x)
 \end{pmatrix}^{*}
   {\mathcal J}_{m-1}^{(m-1)}
  {\mathcal W}_{2m}^{(0)}
    {\mathcal J}_{m-1}^{(m-1)}
 \begin{pmatrix}
 F_0^{(2m-2)}(y)\\
 F_1^{(2m-1)}(y)
 \end{pmatrix}
 \label{F2mCC}
 \end{align}
 and
\begin{align}
 {\mathcal F}_{2m+1}(ix,iy)=&2
 \begin{pmatrix}
 F_0^{(2m)}(\bar x)\\
 F_1^{(2m-1)}(\bar x)
 \end{pmatrix}^{*}
   {\mathcal J}_{m-1}^{(m)}
  {\mathcal W}_{2m+1}^{(0)}
   {\mathcal J}_{m-1}^{(m)}
   \begin{pmatrix}
 F_0^{(2m)}(y)\\
 F_1^{(2m-1)}(y)
 \end{pmatrix}.\label{F2mp1CC}
 \end{align}
\end{corollary}
\begin{proof}
Equalities  (\ref{F2mCC})  and (\ref{F2mp1CC})  readily follow from  (\ref{F2mC})  and (\ref{F2mp1C}), respectively.
\end{proof}
%
Let
\begin{align*}
 V_{[0,2m-1]}(x):=&\left(I_q,\bar xI_q,\bar x^2I_q,\ldots, \bar x^{2m-2}I_q,\bar x^{2m-1}I_q\right)^\mathsf{*},\\
 V_{[0,2m]}(x):=&\left(I_q,\bar xI_q,\bar x^2I_q,\ldots,\bar x^{2m-1}I_q,\bar x^{2m}I_q\right)^\mathsf{*}.
\end{align*}
Furthermore, let ${\mathcal T}_{2m}$ (resp. ${\mathcal T}_{2m+1}$) be elementary block matrices such that 
\begin{equation}\label{TT2m}
 \begin{pmatrix}
 F_0^{(2m-2)}(x)\\
 F_1^{(2m-1)}(x)
 \end{pmatrix}=
 {\mathcal J}_{m-1}^{(m-1)}
 {\mathcal T}_{2m} V_{[0,2m-1]}(x)
 \end{equation}
and 
\begin{equation} \label{TT2mp1}
\begin{pmatrix}
 F_0^{(2m)}(x)\\
 F_1^{(2m-1)}(x)
 \end{pmatrix}=
 {\mathcal J}_{m-1}^{(m)}
 {\mathcal T}_{2m+1} V_{[0,2m]}(x).
 \end{equation}
The following remark readily follows.
\begin{remark}\label{rem5.2}
Using (\ref{TT2m}) and (\ref{TT2mp1}), equalities (\ref{F2mCC}) and (\ref{F2mp1CC}) can be written as follows:
\begin{align}
 {\mathcal F}_{2m}(ix,iy)=&2V_{[0,2m-1]}^*(\bar x){\mathcal T}_{2m}^*
  {\mathcal J}_{m-1}^{(m-1)}
  {\mathcal W}_{2m}^{(0)}
   {\mathcal J}_{m-1}^{(m-1)}
     {\mathcal T}_{2m} V_{[0,2m-1]}(y), \label{VV2m}
     \\
     {\mathcal F}_{2m+1}(ix,iy)=&2
     V_{[0,2m]}^*(\bar x){\mathcal T}_{2m+1}^*
      {\mathcal J}_{m-1}^{(m)}
  {\mathcal W}_{2m+1}^{(0)}
   {\mathcal J}_{m-1}^{(m)}
     {\mathcal T}_{2m+1} V_{[0,2m]}(y).\label{VV2mp1}
 \end{align}
\end{remark}

The positivity of the matrices
${\mathcal W}_{2m}^{(0)}$
and
${\mathcal W}_{2m+1}^{(0)}$
is consistent with the spectral separation described in
Remark~\ref{rem:spectral_symmetries}. Indeed, the substitution
$z=i\lambda$
transforms separation with respect to the imaginary axis into a Hermitian
positivity structure for the associated Bezoutian kernel.

We now state the main result of the present work. The key step is to express
the matrix
\[
\frac{1}{i}B_{L_1^*,L^*}(L,L_1)
\]
in terms of the Bezoutian forms obtained above and to prove its positive
definiteness. This allows us to establish the Hurwitz property of HTM
polynomials.

\begin{theorem}\label{th.3aa}
Each HTM polynomial is a Hurwitz matrix polynomial.
\end{theorem}

\begin{proof} We use Corollary \ref{cor1.19} with
\begin{equation}L(\lambda)=\fop_n(i\lambda) \quad  \mbox{and} \quad L_1(\lambda)=\fop_n(-i\lambda) \label{LL22}\end{equation}
and Equalities (\ref{BM1M}) , (\ref{GML1}), (\ref{fn41}), (\ref{hgn1}), (\ref{hgn2}) and (\ref{FGn}). For $n=2m$, by using (\ref{VV2m}),
we have
\begin{equation}\label{B2m}
{\displaystyle \frac{1}{i}B_{L_1^*,L^*}(L,L_1)}=2{\mathcal T}_{2m}^* {\mathcal J}_{m-1}^{(m-1)}
  {\mathcal W}_{2m}^{(0)}
   {\mathcal J}_{m-1}^{(m-1)}
     {\mathcal T}_{2m}.
\end{equation}
By employing (\ref{VV2mp1}) for $n=2m+1$, we have
\begin{equation}\label{B2mp1}
{\displaystyle \frac{1}{i}B_{L_1^*,L^*}(L,L_1)}=2{\mathcal
T}_{2m+1}^*
 {\mathcal J}_{m-1}^{(m)}
  {\mathcal W}_{2m+1}^{(0)}
   {\mathcal J}_{m-1}^{(m)}
     {\mathcal T}_{2m+1}.
\end{equation}     
     The right-hand sides of (\ref{B2m}) and (\ref{B2mp1}) are positive definite matrices because
     \begin{equation}\label{HHmm}
      \begin{pmatrix}
 H_{2,m-1}&0\\
0& H_{1,m-1}
 \end{pmatrix}
\quad
\mbox{and}
\quad
\begin{pmatrix}
 H_{1,m}&0\\
0& H_{2,m-1}
 \end{pmatrix}
\end{equation}
 are positive definite, and both  ${\mathcal J}_{m-1}^{(m-1)}
     {\mathcal T}_{2m}$ and  ${\mathcal J}_{m-1}^{(m)}
     {\mathcal T}_{2m+1}$ are invertible matrices.
     
     By Corollary~\ref{cor1.19}, the spectra of
\[
L(\lambda)=\fop_n(i\lambda)
\]
belong to the upper half-plane of~$\mathbb C$.
Consequently, the spectra of $\fop_n(z)$ lie in the open left half-plane of~$\mathbb C$.
Hence, the HTM polynomials $\fop_n$ are Hurwitz matrix polynomials.
\end{proof}

\begin{remark}\label{rem5.4}
By employing \eqref{inertia}, note that
\[
\gamma_{-}(\fop_n(i\lambda))
=
\gamma_{0}(\fop_n(i\lambda))
=
0,
\]
and therefore
\[
\gamma(\fop_n(i\lambda))
=
\gamma_{+}(\fop_n(i\lambda)).
\]
This is equivalent to the fact that all zeros of $\det \fop_n(\lambda)$ lie in the open left half-plane.
\end{remark}
The following remark illustrates the previous observation in the simplest nontrivial case.
\begin{remark}\label{rem5.4A}
A $q\times q$ HTM polynomial $\fop_2$ of degree $2$ is a Hurwitz matrix polynomial; it admits the representation
\[
\fop_2(z)=I_q z^2+s_0z+s_0^{-1}s_1.
\]
\end{remark}
The next remark concerns the matrix measures associated with HTM polynomials and their connection with the corresponding orthogonal matrix polynomials.
\begin{remark}
A $q\times q$ matrix $\widetilde \sigma(\cdot)$ defined on $[0,\infty)$ is called monotonically nondecreasing if $\widetilde\sigma(t)$ is self-adjoint for all $t\in[0,\infty)$
and $\widetilde\sigma(t_2)-\widetilde\sigma(t_1)\ge 0$ for $t_1\le t_2$.

For $n=2m+1$ (resp. $n=2m$), 
the points of discontinuity of the distribution function $\widetilde \sigma_{2m,\max}$ (resp. $\widetilde \sigma_{2m,\min}$)
that correspond to the measure $\sigma_{2m,\max}$ (resp. $\sigma_{2m,\min}$)
are given by the roots of the polynomial $\det(P_{1,m}^*(-\bar z))$ (resp. $\det(z P_{2,m}^*(-\bar z))$); see (\ref{GGG1}) and (\ref{GGG2}).

On the other hand,
the polynomials (\ref{h001}) and (\ref{g001}), which are constructed using the moments $s_j$,
can be obtained from an absolutely continuous distribution on $[0,\infty)$ that corresponds to a positive measure on $[0,\infty)$, as shown in the next example.
\end{remark}
The following example illustrates the previous construction in the case of a fifth-degree HTM polynomial.
For the matrix polynomial $\fop_5$, we compute the Bezoutian forms ${\mathcal G}_5^{(1)}$ and ${\mathcal G}_5^{(2)}$, the symmetrizers $S(\gop_5)$ and $S(\widetilde \gop_5)$, the matrices $\widetilde J_1$, $\widetilde J_2$, ${\mathcal J}_1^{(2)}$, as well as the matrices ${\mathcal W}_5^{(0)}$ and ${\mathcal T}_5$. Using these matrices, we obtain an explicit representation of
\[
\frac{1}{i}B_{L_1^*,L^*}(L,L_1),
\]
whose positive definiteness ensures that $\fop_5$ is a Hurwitz matrix polynomial.
\begin{example}
Let
\begin{equation}\label{eq001dd}
\fop_5(z)=z^5 I_2+ A_1 z^4+A_2 z^3+A_3 z^2+A_4 z+A_5,
\end{equation}
where
$A_1=
\begin{pmatrix}
 2 & -i \\
 i & 1
\end{pmatrix},
A_2=
\begin{pmatrix}
 \frac{486}{37} & 0 \\
 -\frac{264}{37}i & 6
\end{pmatrix},
A_3=
\begin{pmatrix}
  \frac{560}{37} & -5i \\
 5i & 5
\end{pmatrix},$\\
$A_4=
\begin{pmatrix}
 \frac{1158}{37} & 0 \\
 -\frac{936}{37}i & 6 \\
\end{pmatrix}$ and
$A_5=
\begin{pmatrix}
 \frac{292}{37} & -2 i \\
 2 i & 2 \\
\end{pmatrix}.
$
\\
The matrix polynomials $\hop_5$ and $\gop_5$ have the form 
\begin{equation}\label{hhgg5}
\hop_5(x)=A_1 x^2+A_3 x+ A_5, \qquad \mbox{and}\qquad \gop_5(x)=I_2 x^2+A_2 x+ A_4
\end{equation}
Using (\ref{hhgg5}) together with the Laurent expansions \eqref{hup09} and \eqref{hup10}, we obtain the first four moments:
\begin{align}
s_0=&\begin{pmatrix}
 2 & -i \\
 i &1
\end{pmatrix},\quad 
s_1=\begin{pmatrix}
 4 & -i \\
 i &1
\end{pmatrix},
s_2=\begin{pmatrix}
 16 & -2i \\
 2i &2
\end{pmatrix},\label{eqss1}
\\
s_3=&\begin{pmatrix}
 96 & -6i \\
 6i &6
\end{pmatrix},\quad 
s_4=\begin{pmatrix}
 768 & -24i \\
 24i &24
\end{pmatrix}.\label{eqss2}
\end{align}
The Bezoutian forms ${\mathcal G}_5^{(1)}$ and ${\mathcal G}_5^{(2)}$ defined in \eqref{hgn1} and \eqref{hgn2} are given by
\begin{align*}
&{\mathcal G}_5^{(1)}(x,y)=\left(
\begin{array}{cc}
g_{11}(x,y)  &-i\,g_{12}(x,y) \\
i\,g_{12}(x,y) & g_{22}(x,y)
\end{array}
\right)
\end{align*}
where
\begin{align*}
g_{11}=&\frac{2}{1369}
\left(
37 x^4 \left(37 y^4+280 y^2+146\right)+4 x^2 \left(2590 y^4+22883 y^2+15297\right)\right.\\
&\left.+5402 y^4+61188 y^2+134436\right)\\
g_{12}=&\left(\left(y^4+5 y^2+2\right) x^4+\left(5 y^4+26
   y^2+12\right) x^2+2 \left(y^4+6 y^2+6\right)\right)\\
g_{22}=& \left(y^4+5 y^2+2\right) x^4+\left(5 y^4+26 y^2+12\right) x^2+2 \left(y^4+6 y^2+6\right) 
\end{align*}
and 
\begin{align*}
&{\mathcal G}_5^{(2)}(x,y)\\
&=\left(
\begin{array}{cc}
 -\frac{4 x y \left(37 \left(37 y^2+272\right) x^2+10064 y^2+88236\right)}{1369} & i x y \left(\left(y^2+4\right) x^2+4 y^2+18\right) \\
 -i x y \left(\left(y^2+4\right) x^2+4 y^2+18\right) & -x y \left(\left(y^2+4\right) x^2+4 y^2+18\right) \\
\end{array}
\right)
\end{align*}
The matrices $S(\gop_5)$ and $S(\widetilde \gop_5)$, defined in \eqref{Sqq}, have the form
\begin{equation}\label{SS551}
S(\gop_5)=\begin{pmatrix}A_4&A_2&I_2\\
A_2&I_2&0_2\\
I_2&0_2&0_2
\end{pmatrix}=
\left(
\begin{array}{cccccc}
 \frac{1158}{37} & 0 & \frac{486}{37} & 0 & 1 & 0 \\
 -\frac{936 i}{37} & 6 & -\frac{264 i}{37} & 6 & 0 & 1 \\
 \frac{486}{37} & 0 & 1 & 0 & 0 & 0 \\
 -\frac{264 i}{37} & 6 & 0 & 1 & 0 & 0 \\
 1 & 0 & 0 & 0 & 0 & 0 \\
 0 & 1 & 0 & 0 & 0 & 0 \\
\end{array}
\right)
\end{equation}
and
\begin{equation}\label{SS552}
S(\widetilde \gop_5)=\begin{pmatrix}A_2&I_2\\
I_2&0_2
\end{pmatrix}=
\left(
\begin{array}{cccccc}
 \frac{1158}{37} & 0 & \frac{486}{37} & 0 & 1 & 0 \\
 -\frac{936 i}{37} & 6 & -\frac{264 i}{37} & 6 & 0 & 1 \\
 \frac{486}{37} & 0 & 1 & 0 & 0 & 0 \\
 -\frac{264 i}{37} & 6 & 0 & 1 & 0 & 0 \\
 1 & 0 & 0 & 0 & 0 & 0 \\
 0 & 1 & 0 & 0 & 0 & 0 \\
\end{array}
\right).
\end{equation}
The matrix polynomial $\tilde \gop_5$ has the form
$$
\tilde \gop_5(x)= I_2 x+A_2.
$$
Furthermore, the matrices $F_0^{(4)}$ and $F_1^{(3)}$, defined in \eqref{F02k} and \eqref{F12km1}, have the form
\begin{equation}\label{FF43}
F_0^{(4)}(x)=\begin{pmatrix}I_2&\bar x^2 I_2&\bar x^4 I_2\end{pmatrix}^*\quad \mbox{and}\quad
F_1^{(3)}(x)=\begin{pmatrix}I_2&\bar x^2 I_2\end{pmatrix}^*.
\end{equation}
The matrices $\widetilde J_2$ and $\widetilde J_1$, introduced in \eqref{Jqq}, take the form
\begin{equation}\label{JJ21}
\widetilde J_2=\begin{pmatrix}
I_2& 0_2&0_2\\
0_2& -I_2&0_2\\
0_2& 0_2&1_2
\end{pmatrix},\quad \mbox{and}\quad
\widetilde J_1=\begin{pmatrix}
I_2& 0_2\\
0_2& -I_2
\end{pmatrix}.
\end{equation}
Using \eqref{eqss1} and \eqref{eqss2}, we construct the matrices $H_{1,2}$ and $H_{2,1}$.
One readily verifies that these block matrices are positive definite.
From \eqref{SS551}, \eqref{SS552}, \eqref{FF43}, and \eqref{JJ21},
we obtain the right-hand sides of \eqref{GGa2mp1} and \eqref{GGa2mp2}, which allow us to conclude that
${\mathcal G}_5^{(1)}$ and $-{\mathcal G}_5^{(2)}$ are positive quadratic forms.
\\
The matrix $\frac{1}{i}B_{L_1^*,L^*}(L,L_1)$ in \eqref{B2mp1} can be written as the product of the matrices
$${\mathcal J}_1^{(2)}=\begin{pmatrix}
\widetilde J_2&0\\
0& \widetilde J_1
\end{pmatrix},$$
\begin{align*}
{\mathcal W}_{5}^{(0)}=&
  \begin{pmatrix}
 S(\gop_{5}^*)&0\\
0& S(\widetilde \gop_{5}^*)
 \end{pmatrix}
  {\mathcal J}_{1}^{(2)}
 \begin{pmatrix}
 H_{1,2}&0\\
0&  H_{2,1}
 \end{pmatrix}
   {\mathcal J}_{1}^{(2)}
  \begin{pmatrix}
 S(\gop_{5})&0\\
0& S(\widetilde \gop_{5})
 \end{pmatrix}\label{WW2A}
\end{align*}
and the matrix ${\mathcal T}_5$, which appears in \eqref{VV2mp1}, has the form
$${\mathcal T}_5=\begin{pmatrix}
I_2&0_2&0_2&0_2&0_2&0_2\\
0_2&0_2&I_2&0_2&0_2&0_2\\
0_2&0_2&0_2&0_2&0_2&I_2\\
0_2&I_2&0_2&0_2&0_2&0_2\\
0_2&0_2&0_2&0_2&I_2&0_2
\end{pmatrix}.$$
From these matrices, together with \eqref{VV2mp1}, \eqref{B2mp1}, Corollary~\ref{cor1.19}, and the second matrix in \eqref{HHmm},
it follows that the matrix polynomial $\fop_5$ is a Hurwitz matrix polynomial.
\end{example}
The following remark illustrates how HTM polynomials can be generated from matrix-valued measures on $[0,\infty)$.
\begin{remark} 
\label{remx01}
The polynomial $\fop_5$ in \eqref{eq001dd} can be constructed starting from the matrix-valued distribution
\[
\widetilde \sigma(x)=
\begin{pmatrix}
 4-2 e^{-x/2} & i e^{-x} \\
 -i e^{-x} & 2-e^{-x}
\end{pmatrix},
\]
defined on $[0,\infty)$.
\\
More generally, the construction of HTM polynomials can be initiated from a nonnegative Hermitian matrix-valued measure (or a distribution of bounded variation) on $[0,\infty)$ such that the associated block Hankel matrices $H_{1,j}$ and $H_{2,j}$ are positive definite.
\\
From this distribution, we compute the moments $s_j$ as in \eqref{eqss1} and \eqref{eqss2}. Next, we apply Definition~\ref{def2p}. Finally, the second equality in \eqref{hupfqp} is used to obtain \eqref{eq001dd}.
\end{remark}
Next, we recall the characterization of HTM polynomials in terms of block Hankel matrices; see \cite[Theorem 7.10]{abH}.
\begin{theorem} \label{th2.2A}
Let $n$ be a natural number with $n\ge 2$, and let $\fop_n$ be a matrix-valued
polynomial. Let $\hop_n$ (resp. $\gop_n$) denote the $\hop_n$-part
(resp. $\gop_n$-part) of $\fop_n$.
Furthermore, let $n=2m$ (resp. $n=2m+1$) with $m\ge 1$, and let
$H_{1,m-1}$ and $H_{2,m-1}$ (resp. $H_{1,m}$ and $H_{2,m-1}$)
be constructed from the Hermitian Markov parameters of the polynomial $\fop_n$.
Then the matrix-valued polynomial $\fop_n$ is an HTM polynomial if and only if
the block Hankel matrices $H_{1,m-1}$ and $H_{2,m-1}$
(resp. $H_{1,m}$ and $H_{2,m-1}$) are positive definite.
\end{theorem}
We now pass from the Hankel-matrix characterization to a Bezoutian characterization.
In particular, for an HTM polynomial, the Bezoutians associated with the decomposition
of $\fop_n$ into its $\hop_n$- and $\gop_n$-parts have definite signs.
\begin{proposition}\label{prop5-4}
Let $\fop_n$ be an HTM polynomial. Moreover, let $\hop_n$ and $\gop_n$ be as in (\ref{hup02}),  (\ref{hn}) and (\ref{gn}), respectively.
Furthermore, let  $B_{x\gop_{n}^*,\hop_{n}^*}(x\gop_{n},\hop_{n})$ and
$B_{\gop_{n}^*,\hop_{n}^*}(\gop_{n},\hop_{n})$ be Bezoutians associated with the quadruple $(x\gop_{n}^*,\hop_{n},\hop_{n}^*,x\gop_{n})$
and $(\gop_{n}^*,\hop_{n},\hop_{n}^*,\gop_{n})$, respectively. Thus,
\\
a) the matrix
$B_{x\gop_{n}^*,\hop_{n}^*}(x\gop_{n},\hop_{n})$ is  positive definite, and\\
b) the matrix $B_{\gop_{n}^*,\hop_{n}^*}(\gop_{n},\hop_{n})$ is  negative definite.
\end{proposition}
\begin{proof}
Let $\widetilde J_m$, $\widetilde \gop_{2m+1}$,
 $S(\hop_{2m})$, $S(\gop_{2m+1})$ and $S(\widetilde \gop_{2m+1})$ be as in (\ref{Jqq}), (\ref{tg2m1}) and  (\ref{Sqq}).
 For $n=2m$, by  using (\ref{hgn1}), (\ref{hgn2}), (\ref{GGa2m1}) and (\ref{GGa2m2}), 
  we have
\begin{align*}
B_{x \gop_{2m}^*,\hop_{2m}^*}(x\gop_{2m},\hop_{2m})=&S(\hop_{2m}^*) \widetilde J_{m-1} H_{2,m-1}\widetilde J_{m-1}S(\hop_{2m}),\\
B_{\gop_{2m}^*,\hop_{2m}^*}(\gop_{2m},\hop_{2m})=&-S(\hop_{2m}^*)\widetilde
J_{m-1} H_{1,m-1}\widetilde J_{m-1} S(\hop_{2m}).
\end{align*}
By (\ref{hgn1}), (\ref{hgn2}), (\ref{GGa2mp1}) and (\ref{GGa2mp2}), for $n=2m+1$,  we obtain
\begin{align*}
B_{x\gop_{2m+1}^*,\hop_{2m+1}^*}(x\gop_{2m+1},\hop_{2m+1})=&S(\gop_{2m+1}^*)\widetilde J_{m} H_{1,m}\widetilde J_{m} S(\gop_{2m+1})\\
B_{\gop_{2m+1}^*,\hop_{2m+1}^*}(\gop_{2m+1},\hop_{2m+1})=&-S(\widetilde
\gop_{2m+1}^*)\widetilde J_{m-1} H_{2,m-1}\widetilde J_{m-1}
S(\widetilde\gop_{2m+1}).
\end{align*}
From Theorem \ref{th2.2A}, see that assertion a) and b) are valid.
 \end{proof}
\begin{remark}\label{rem5.6}
In  \cite[Theorem 2]{lan}, it is proven that a scalar polynomial $f_n$ is Hurwitz if and only if
$B_{x g_{n}^*,h_{n}^*}(xg_{n},h_{n})$ is  negative definite and
 $B_{g_{n}^*,h_{n}^*}(g_{n},h_{n})$ is  positive definite.
 Note that the Markov parameters related to \cite[Theorem 2]{lan} have the opposite sign compared with those associated with parts (a) and (b) of
 Proposition \ref{prop5-4}.
\end{remark}


\section{Completing a matrix polynomial to an HTM polynomial}\label{sec3aa}
\medskip
As indicated in \cite{abcomment}, not every Hurwitz matrix polynomial is an HTM polynomial. In this section, we propose a procedure to complete a given matrix polynomial into an HTM polynomial.

In \cite{abcomment}, the polynomial
\begin{equation}\label{eq001a}
P_3(z)=I_2 z^3 + A_1 z^2 + A_2 z + A_3
\end{equation}
with
\begin{equation}\label{AA13}
A_1=\begin{pmatrix}
\frac{2187}{109} & \frac{1206}{109} \\
0 & 9
\end{pmatrix},\quad
A_2=\begin{pmatrix}
\frac{10662}{109} & \frac{8700}{109} \\
0 & 18
\end{pmatrix},
\quad
A_3=\begin{pmatrix}
\frac{11178}{109} & \frac{10524}{109} \\
0 & 6
\end{pmatrix}
\end{equation}
was suggested as a polynomial that is Hurwitz but not
an HTM polynomial. The rational matrix function
$\frac{1}{z}\left(A_1 z + A_3\right)\left(I_2 z + A_2\right)^{-1}$
does not satisfy property~(\ref{sfn}). Moreover, as indicated in \cite{abcomment}, the $2\times 2$ matrix coefficients $c_j$ of the expansion
\[
\frac{1}{z}\left(A_1 z + A_3\right)\left(I_2 z + A_2\right)^{-1}
= -\frac{c_0}{z}-\frac{c_1}{z^2}-\frac{c_2}{z^3}-\cdots
\]
are not Hermitian matrices. Consequently, the polynomial $P_3$ cannot be investigated in the framework of HTM polynomials.

The following result allows one to determine whether a matrix polynomial $P_n$ that is not an HTM polynomial is in fact a Hurwitz polynomial by ``completing,'' mainly by perturbing $P_n$ with another polynomial $Q_{n-1}$. This is done by constructing a new HTM polynomial of degree $2n$.
\begin{theorem}\label{complet}
Let $P_n(z)$ be $q\times q$ a monic polynomial of degree $n$. If
there exists a  $q\times q$ polynomial  $Q_{n-1}(z)$ of degree $n-1$
such that
\begin{equation}\label{f2npq}
\fop_{2n}(z)=P_n(z^2)+zQ_{n-1}(z^2)
\end{equation}
 is an HTM polynomial, then the polynomial $P_n(z)$ is a Hurwitz polynomial.
\end{theorem}
\begin{proof}
From (\ref{hupfqp}) and (\ref{h001}), there are polynomials
$P_{1,n}$ y $Q_{1,n}$ such that $P_n(z)=(-1)^nP_{1,n}^*(-\bar z)$
$Q_{n-1}(z)=(-1)^{n+1}Q_{1,n}^*(-\bar z)$ and $P_{1,n}(z)$ is an
orthogonal polynomial on $[0,+\infty)$ with respect to a certain
positive measure $\sigma$ corresponding to the Stieltjes transform
$
-\frac{Q_{1,n}^*(\bar z)}{P_{1,n}^*(\bar z)}$. By \cite[Theorem
1]{dyu2018}, the roots of the polynomial $\det P_{1,n}(z)$ lie in
$[0,+\infty)$.
Consequently, the roots of the polynomial
\[
\det P_n(z)=\det\bigl((-1)^n P_{1,n}^*(-\bar z)\bigr)
\]
lie in the open left half-plane. Therefore, $\det P_n(z)$ is a Hurwitz polynomial.
\end{proof}

\begin{remark}\label{rem4.5}
Let $P_n$ be a monic matrix polynomial and let $P_n^\intercal$ denote its transpose.
Since
\[
\det P_n=\det P_n^\intercal,
\]
Theorem~\ref{complet} is also valid for
$P_n^\intercal$ and $Q_{n-1}^\intercal$.
\end{remark}
For the construction of the polynomial $Q_{n-1}$ appearing in Theorem~\ref{complet}, we propose the following algorithm. 
The procedure is based on the representation of the associated rational matrix function in terms of orthogonal polynomials and the corresponding Markov parameters, which allows one to recover $Q_{n-1}$ from a given polynomial $P_n$.
\\
{\bf Algorithm for finding $Q_{n-1}$ from $P_{n}$}
\\
Step 1. Let  $P_{1,n}(z)=(-1)^n P_n^*(-\bar z)$ (resp. $P_{1,n}(z)=(-1)^n P_n^{\intercal^*}(-\bar z)$). By using (\ref{qp1}), compute $Q_{1,n}(z)$
in terms of the coefficients of
 $P_n(z)$ (resp. $P_n^\intercal(z)$)
and the moments or Markov parameters $s_j$,
associated with some positive measure  $\sigma$
  on $[0,+\infty)$.\\
Step 2. Verify that $\frac{P_{1,n}^*(\bar z)}{Q_{1,n}^*(\bar z)}$
 satisfies Equality (\ref{sfn}).  If Equality (\ref{sfn}) is not satisfied, we cannot apply this algorithm.
\\
Step 3. By using the rational matrix function $ \frac{P_{1,n}^*(\bar z)}{Q_{1,n}^*(\bar z)}$,
calculate the remaining Markov parameters from the set  $(s_j)_{j=0}^{2n-1}$.
\\
Step 4. If sequence $(s_j)_{j=0}^{2n-1}$ is a positive definite Stieltjes sequence, meaning that the corresponding
 Hankel matrices $H_{1,n-1}$ and $H_{2,n-1}$ are positive definite matrices, then $Q_{n-1}(z)=(-1)^{n+1}Q_{1,n}^*(-\bar z)$.
\begin{example}
We apply the completion procedure to the matrix polynomial (\ref{eq001a}).

For $P_3^\intercal(z)$, we use the above algorithm together with Remark \ref{rem4.5}.   
With moments
$s_0=\begin{pmatrix}
 2 b & -b \\
 -b & b
\end{pmatrix}
$,
$s_1=\begin{pmatrix}
 4 b & -b \\
 -b & b
\end{pmatrix}
$
 and
$s_2=\begin{pmatrix}
 16 b & -2b \\
 -2b & 2b
\end{pmatrix}
$,
we calculate $Q_{1,3}(z)$; therefore,  we have
$$
Q_{1,3}(z)=s_0 z^2+ z(s_1-A_1^{\intercal^*} s_0 )+s_2-A_1^{\intercal^*} s_1+A_2^{\intercal^*} s_0 .
$$
From the Laurent expansion of $-\frac{P_{1,3}^*(\bar z)}{Q_{1,3}^*(\bar z)}$, we calculate the moments
$$s_3=\begin{pmatrix}
 96 b & -6b \\
 -6b & 6b
\end{pmatrix}
,
s_4=\begin{pmatrix}
 768 b & -24b \\
 -24b & 24b
\end{pmatrix}
,
s_5=\begin{pmatrix}
 7680 b & -120b \\
 -120b & 120b
\end{pmatrix}
.$$
 The block Hankel  matrices $H_{1,2}$ and $H_{2,2}$ are positive
definite matrices for $b>0$. From (\ref{h001}),  polynomial $Q_2(z)$
constructed from $P_3$ has the form $Q_2(z)=Q_{1,3}^{*}(-\bar z)$.
Thus, the resulting HTM polynomial $\fop_6$ can be represented as in
(\ref{f2npq}):
\begin{equation}\label{f6exp}
\fop_6(z)=z^6I_2+B_1z^5+A_1z^4+B_2z^4+A_2z^2+B_3z+A_3.
\end{equation}
The coefficients $A_1$, $A_2$ and $A_3$ are given by (\ref{AA13}), and the coefficients $B_1$, $B_2$ and $B_3$ are
$$
B_1=\begin{pmatrix}
 2b & -b \\
 -b & b
\end{pmatrix},
B_2=\begin{pmatrix}
 \frac{2732 b}{109} & -8 b \\
 -8 b & 8 b \\
\end{pmatrix},
B_3=\begin{pmatrix}
\frac{6826 b}{109} & -11 b \\
 -11 b & 11 b
\end{pmatrix}.
 $$
By Theorem \ref{th.3aa}, matrix polynomial (\ref{f6exp}) is a Hurwitz
matrix polynomial for $b>0$. Finally, by Theorem \ref{complet},
polynomial $P_3$ as in (\ref{eq001a}) is a Hurwitz matrix
polynomial.
\end{example}
\vskip2mm
In the scalar version, one can construct a sequence of scalar Hurwitz polynomials  $(f_j)_{j=1}^\infty$
using the sequence of scalar orthogonal polynomials on $[0,+\infty)$   $(p_{1,j})_{j=1}^\infty$ and $(p_{2,j})_{j=1}^\infty$
as in Definition \ref{def2p}; see \cite{ab2020A}.

In the following lemma, we propose constructing a family of scalar Hurwitz polynomials involving the polynomials
$p_{1,j}$, and $p_{2,j}$ and their derivatives.
\begin{lemma}\label{lem4.7} For $q=1$, let
 $(p_{1,m}(z))$ (resp. $(p_{2,m}(z))$) be as in Definition \ref{def2p}, and let
 $\widetilde p_{1,m-1}(z):=\frac{d}{dz}p_{1,m}(z)$ (resp. $\widetilde p_{2,m-1}(z):=\frac{d}{dz}p_{2,m}(z)$). Thus,
for $m\geq2$, the polynomials
 \begin{equation}\label{f2m001}
  f_{2m}(z)=(-1)^m(p_{1,m}(-z^2)-z \widetilde p_{1,m-1}(-z^2))
  \end{equation}
  and
   \begin{equation}\label{f2m002}
  f_{2m+1}(z)=(-1)^m(\widetilde p_{2,m}(-z^2)+z p_{2,m}(-z^2)).
  \end{equation}
  are scalar Hurwitz polynomials.
\end{lemma}
\begin{proof}
The proof of (\ref{f2m001}) and (\ref{f2m002}) readily follows from the interlacing properties of the roots of the orthogonal polynomials 
\( p_{k,m}(z) \) and \( \widetilde{p}_{k,m-1}(z) \) for \( k = 1, 2 \), and from the fact that \( a_0 a_1 > 0 \), where \( a_0 \) and \( a_1 \) are the leading coefficients of the polynomial
\( f_n(z) = a_{0} z^n + a_{1} z^{n-1} + \cdots + a_{n} \);
see the Hermite--Biehler theorem \cite{Biehler}, \cite{Hermite}.
\end{proof}

\begin{remark}\label{rem5.14}
From (\ref{eq001dd}) and (\ref{f6exp}), we see that the determinants of the coefficients of these polynomials are positive.
Based on this observation, we propose the conjecture below.
\\
{\bf Conjecture}.
The determinants of the coefficients $A_j$ of an HTM polynomial $\fop_n$ are positive.
\end{remark}

\section{Additional remarks}
\label{sec06}

In this section, to clarify the importance of the present work,
we recall several relevant results from \cite{zhan1},
paying particular attention to certain steps in the proofs
that were left implicit.
We emphasize that the authors of \cite{zhan1} (2021) were the first
to propose that the HTM polynomials possess the Hurwitz property,
while the concept of HTM polynomials itself was introduced earlier
by the author of the present work in \cite{abH} (2015).

\begin{lemma}\cite[Lemma 3.4]{zhan1}\label{lemZ1}
Let $N_R(z)$, $D_R(z)$, $N_L(z)$ and $D_L(z)\in\mathbb{C}[z]^{p\times p}$ 
be regular such that $\deg D_R > \deg N_R$, $\deg D_L > \deg N_L$ and,
for all large enough $z\in\mathbb{C}$,
\[
(D_L(z))^{-1}N_L(z)=N_R(z)(D_R(z))^{-1}
=\sum_{k=0}^{\infty} z^{-(k+1)}s_k .
\]

If $D_L(z)$ and $D_R(z)$ are written in the form
\[
D_L(z)=\sum_{k=0}^{m_L} D_{L,m_L-k}z^k
\qquad \text{and} \qquad
D_R(z)=\sum_{k=0}^{m_R} D_{R,m_R-k}z^k ,
\]
where $D_{L,k}\in\mathbb{C}^{p\times p}$ for $k=0,\ldots,m_L$ and
$D_{R,k}\in\mathbb{C}^{p\times p}$ for $k=0,\ldots,m_R$, then
\begin{align}
\mathbf{B}_{D_L,N_L}(D_R,N_R)
=&
\begin{pmatrix}
D_{L,m_L-1} & \cdots & D_{L,0} \\
\vdots & \ddots & \vdots \\
D_{L,0} & \cdots &
\end{pmatrix}
\begin{pmatrix}
s_0 & \cdots & s_{m_R-1} \\
\vdots & \ddots & \vdots \\
s_{m_L-1} & \cdots & s_{m_R+m_L-2}
\end{pmatrix}\nonumber\\
&\cdot
\begin{pmatrix}
D_{R,m_R-1} & \cdots & D_{R,0} \\
\vdots & \ddots & \vdots \\
D_{R,0} & \cdots &
\end{pmatrix}. \label{BB01}
\end{align}
\end{lemma}

Let  $L_1$ and $L$ be as in (\ref{LL22}).  In \cite{zhan1}, the Theorem 3.5 is regarding the inertia representation of matrix polynomials $\fop_n$.  
For $n=2m$,
in the proof of \cite[Theorem 3.5]{zhan1} the \cite[Equality (3.6)]{zhan1} in terms of our notations is the following
 \begin{equation}
B_{L_1^*,L^*}(L,L_1)=2iB_{\hop_{2m}^*(-\bar x^2),x \gop_{2m}^*(-\bar x^2)}(\hop_{2m}(-x^2),-x \gop_{2m}(-x^2)). \label{LLBB}
 \end{equation}
The authors of \cite{zhan1} use $2i\,B_{F_e^{\vee},\,-F_0^{\vee}}
\bigl(\widehat{F}_e,\,-\widehat{F}_0\bigr)$ in place of the right-hand side of \eqref{LLBB}.
They claim that Lemma~\ref{lemZ1} implies that the matrix $2i\,B_{F_e^{\vee},\,-F_0^{\vee}}
\bigl(\widehat{F}_e,\,-\widehat{F}_0\bigr)$ is congruent to the block Hankel matrix
\begin{equation}\label{eqSS0}
\begin{pmatrix}
s_0 & 0_p & s_1 & \cdots & s_{m-1} & 0_p \\
0_p & s_1 & 0_p & \cdots & 0_p & s_m \\
s_1 & 0_p & s_2 & \cdots & s_m & 0_p \\
\vdots & \vdots & \vdots & \ddots & \vdots & \vdots \\
s_{m-1} & 0_p & s_m & \cdots & s_{2m-2} & 0_p \\
0_p & s_m & 0_p & \cdots & 0_p & s_{2m-1}
\end{pmatrix},
\end{equation}
and that, after a suitable reordering of rows and columns in \eqref{eqSS0}, it is further congruent to
\begin{equation}\label{eqHH01}
\begin{pmatrix}
H_{1,m-1} & 0 \\
0 & H_{2,m-1}
\end{pmatrix}.
\end{equation}

\begin{remark}\label{rem5.7A}
By comparing the Hankel matrix of the right-hand side of (\ref{BB01})  and the matrix (\ref{eqSS0})
It is not clear how from  $2i\,B_{F_e^{\vee},\,-F_0^{\vee}}
\bigl(\widehat{F}_e,\,-\widehat{F}_0\bigr)$ one gets (\ref{eqSS0}).
Furthermore, it is not clear how from (\ref{eqSS0}) one obtains (\ref{eqHH01}).\\
In contrast, for $n=2m$, in the present work we derive the matrix
$\left(\begin{smallmatrix}
H_{2,m-1} & 0 \\
0 & H_{1,m-1}
\end{smallmatrix}\right)$
associated with the quadratic form in \eqref{eqHH01}
by employing the decomposition \eqref{FGn} and by computing explicitly
the quadratic forms corresponding to \eqref{hgn1} and \eqref{hgn2},
as shown in \eqref{F2mC}.
Additionally, 
for $n=2m+1$, the proof of \cite[Theorem 3.5]{zhan1} is not given. 
\end{remark}

In Theorem 4.10 of \cite{zhan1}, a ``stability criterion via Markov parameters''  for matrix polynomials
$\fop_n$ in even and odd cases is suggested. We have the following remark regarding this.
\begin{remark}\label{rem5.8A}
For $n=2m$,   Theorem 4.10 of \cite{zhan1}, is related to the HTM polynomials. Other matrix polynomials, for example, the Hurwitz polynomial (\ref{eq001a}), are not covered by \cite[Theorem 4.10]{zhan1}.\\
In this work, we do not discuss the case $n=2m+1$ of Theorem 4.10 of \cite{zhan1}.
\end{remark}


\section{Conclusion}

\medskip

In this work we obtained an explicit representation of the Bezoutian associated with Hurwitz-type matrix polynomials. This representation clarifies the relation between Hurwitz-type matrix polynomials and classical Hurwitz matrix polynomials and provides a constructive tool for verifying Hurwitzness in the matricial setting. Our approach, based on the decomposition of Bezoutians and on the structure of truncated matricial Stieltjes moment problems, offers a natural framework for extending classical stability criteria to matrix polynomials.
In addition, we proposed a method to enlarge the class of Hurwitz-type matrix polynomials by perturbing a polynomial that is not of Hurwitz type so that the new polynomial becomes Hurwitz type.
These results contribute to a better understanding of the structure of
Hurwitz-type matrix polynomials and their relation to Bezoutians,
and they suggest further developments in matrix moment problems,
continued fractions, and stability theory.

\end{document}